\documentclass[10pt]{amsart}
\usepackage[margin=1.15in]{geometry}
\usepackage{amscd,amssymb, amsmath, wasysym}
\usepackage{graphicx}
\usepackage{amsfonts}
\usepackage{mathrsfs}    
\usepackage{amsmath}    
\usepackage{amsthm}     
\usepackage{amscd}      
\usepackage{amssymb}    
\usepackage{eucal}      
\usepackage{latexsym}   
\usepackage{graphicx}   
\usepackage{verbatim}   
\usepackage[all]{xy}     

\usepackage{bookmark}

 \usepackage{hyperref}
 \hypersetup{
     colorlinks=true,
     linkcolor=red,
     filecolor=blue,
     citecolor = blue,      
     urlcolor=cyan,
  }

\makeatletter

\newcounter{thmcounter}

\numberwithin{thmcounter}{section}
\numberwithin{equation}{thmcounter}

\newtheorem{theorem}[thmcounter]{Theorem}
\newtheorem{proposition}[thmcounter]{Proposition}
\newtheorem{lemma}[thmcounter]{Lemma}
\newtheorem{corollary}[thmcounter]{Corollary}

\theoremstyle{definition}
\newtheorem{definition}[thmcounter]{Definition}
\newtheorem{localsetup}[thmcounter]{Local Setup}

\newtheorem{multi-index}[thmcounter]{Multi-index Notation}

\newtheorem{remark}[thmcounter]{Remark}

\newtheoremstyle{claim}{9pt}{3pt}{}{\parindent}{\bf}{.}{1em}{}

\theoremstyle{claim}
\newtheorem{claim}[equation]{Claim}

\newenvironment{namelist}[1]{%
\begin{list}{}
{
\settowidth{\labelwidth}{#1}%
\setlength{\labelsep}{0.3em}%
\setlength{\leftmargin}{\labelwidth}%
\addtolength{\leftmargin}{\labelsep}}}{%
\end{list}}

\newcommand{\nC}{\mathbb{C}}                     

\newcommand{\nG}{\mathbb{G}}                     
\newcommand{\nP}{\mathbb{P}}                     

\newcommand{\nA}{\mathbb{A}}                     
\newcommand{\bT}{\mathbf{T}}  
\newcommand{\bU}{\mathbf{U}}

\newcommand{\bF}{\mathbf{F}} 
\newcommand{\bII}{\mathbf{II}} 
\newcommand{\bIII}{\mathbf{III}} 

\newcommand{\sF}{\mathscr{F}}

\newcommand{\sO}{\mathscr{O}}                    

\newcommand{\ud}{\partial} 

\newcommand{\mf}[1]{\mathfrak{#1}}

\DeclareMathOperator{\codim}{codim}              

\DeclareMathOperator{\Hom}{Hom}                  

\DeclareMathOperator{\im}{im}                    
\DeclareMathOperator{\uIm}{Im}                   

\DeclareMathOperator{\Ker}{Ker}                  

\DeclareMathOperator{\length}{length}            

\DeclareMathOperator{\Sec}{Sec}                  
\DeclareMathOperator{\Sym}{Sym}                  
\DeclareMathOperator{\sHom}{\mathscr{H}om}       
\DeclareMathOperator{\Spec}{Spec}                

\DeclareMathOperator{\Tan}{Tan}                  

\DeclareMathOperator{\rank}{rank}                

\newcommand*{\verylongrightarrow}{\ensuremath{\joinrel\relbar\joinrel\relbar\joinrel\relbar\joinrel\relbar\joinrel\rightarrow}}

\newcounter{rkcounter}             
\makeatletter
\@namedef{subjclassname@2020}{%
	\textup{2020} Mathematics Subject Classification}
\makeatother

\begin{document}

\title[]{On vanishing of fundamental forms of algebraic varieties}
\author{Lawrence Ein}
\address{Department of Mathematics, University Illinois at Chicago, 851 South Morgan St.,
	Chicago, IL 60607, USA}
\email{ein@uic.edu}

\author{Wenbo Niu}
\address{Department of Mathematical Sciences, University of Arkansas, Fayetteville, AR 72701, USA}
\email{wenboniu@uark.edu}

\date{\today}


\begin{abstract} We study fundamental forms of algebraic varieties using the sheaves of principal parts of line bundles and establish a vanishing theorem for any order fundamental forms. We also give connection of fundamental forms with the higher order Gauss map and higher order tangent varieties. 

\end{abstract}

\maketitle


\section{Introduction}

\noindent Throughout we work over the field of complex numbers. Let $X\subseteq \nP^r$ be a {\em quasi-projective} variety of dimension $n\geq 1$. In the influential work \cite{Griffiths:AGandLDG} by Griffiths-Harris, using the method of moving frames one may attach to a generic point $x$ of $X$ a sequence of linear systems of quadrics, cubics, etc, on the projectivization space $\nP(T_x(X))$ of the Zariski tangent space at $x$. These linear systems are called fundamental forms of $X$ in $\nP^r$ and they have deep connections with  the local and global geometry of $X$. We refer the reader to the introduction  of \cite{Griffiths:AGandLDG} for more details. 

Among other things, the degeneration  phenomena from the tangential variety and the secant variety, as well as the Gauss map associated to $X$ are closely related to vanishing of fundamental forms, as studied in \cite{Griffiths:AGandLDG} and in a series work \cite{Landsberg:SecFundFrm,Landsberg:DegSecTan} by Landsberg.  This essentially relates to another fundamental question that how one can determine if a variety is contained in a linear space. To be precise, if $X$ is projective, Fulton-Hansen theorem \cite[5.5]{Lazarsfeld:ConnAppAlgGeo} says that if the tangential variety $\Tan(X)$ (the union of tangent stars) or the secant variety $\Sec(X)$ (the Zariski closure of the union of lines passing through two distinct points of $X$) does not have the expected dimension, then they must be equal and this means that the secant variety is degenerate. Griffiths-Harris' result \cite[6.15]{Griffiths:AGandLDG} further shows that if $\Sec(X)$ degenerates but $\Tan(X)$ does not, then the third fundamental form of $X$ vanishes. Landsberg \cite[10.2]{Landsberg:DegSecTan} pushed this one more step by showing that if $X$ is nonsingular and $\Sec(X)$ degenerates then the third fundamental form vanishes. Note that the vanishing of the third fundemantal form implies that $X$ is contained in its generic second projective tangent space. These results play an important role in the classification of varieties. 

The motivation and the first main purpose of this paper is to reveal the full geometric picture behind the aforementioned classic results and establish a vanishing theorem for arbitrary order fundamental forms. To achieve this and as the second main purpose, instead of using moving frames, we develop an algebraic foundation to the theory of  fundamental forms using sheaves of principal parts, along the line of \cite[Chapter 16]{Grothendieck:EGAIV1967}. Based on our approach, we also show connection between fundamental forms and higher order Gauss maps and higher order tangent varieties (related subjects have been studied previously by various methods, for instance, in  \cite{Piene:HiDualVar}, \cite{Zak:TanSecAlgVar}, \cite{Poi:HiOscuDef}, \cite{Poi:HiGaussMap} and \cite{Rocco:HiGaussMap}). We hope further geometry application based on our approach will come out in near future. 

Let $L=\sO_X(1)$ and  $V=H^0(\nP^r,\sO_{ \nP^{r}}(1))$ be the space of linear forms. For each integer $k\geq 0$, the sheaf of $k$-th order principal parts $P^k(L)$ of $L$ is equipped with a Taylor series map  $$\alpha_k:V\otimes \sO_X\longrightarrow P^k(L).$$ We consider two sheaves $R_k$ and $P_k$  such that $R_k\otimes L$ is the kernel and $P_k\otimes L$ is the image of $\alpha_k$. Over the nonsingular locus of $X$, the kernel of the truncation map $P^k(L)\rightarrow P^{k-1}(L)$ is $S^{k}\Omega^1_X\otimes L$, where $S^k\Omega_X^1$ is the $k$-th symmetric product of the sheaf of differentials. We define the $k$-th fundamental form to be the induced morphism  
$$\bF^*_k:R_{k-1}\longrightarrow S^k\Omega^1_X.$$
If $x\in X$ is a nonsingular point, tensoring with the residue field $k(x)$ one obtains a morphism $$\bF^*_{k,x}:R_{k-1,x}(X)\longrightarrow S^k\Omega^1_x(X)$$
whose image gives a degree $k$ linear system, denoted by $|\bF_{k,x}|$, on the space $\nP(T_x(X))$. Considering $R_{k-1}$ as subsheaf of the conormal sheaf $N^*_X$ of $X$, one can deduce a local formula for the fundamental form $\bF^*_k$ (Corollary \ref{p:04}). The sheaf filtration $N^*_X\supseteq R_2\supseteq R_3\supseteq\ldots$ is strictly decreasing either to zero or to a trivial locally free sheaf after at most $\codim X$ steps (Proposition \ref{p:20}). Thus the $k$-th fundamental forms eventually becomes zero if $k$ is larger than $\codim X$. It turns out that the vanishing of the $k$-th fundamental form implies that the variety $X$ is contained in its generic $(k-1)$-th projective tangent space (the implication is well-known and proved in differential geometry but one can prove it by the method of this paper). 

Geometrically, the sheaf $P_k\otimes L$ gives the $k$-th projective tangent spaces (also called osculating spaces) $\bT^k_x(X)$ of dimension $\mf{t}_k$ at a generic point $x$  (we write $\widehat{\bT}^{k}_x(X)$ for the deprojectivization in $V^*$). For instance, $P_1\otimes L=P^1(L)$ gives projective tangent spaces. It is well known that the second fundamental form is induced by the differential of the Gauss map of $X$. Inspired by Kleiman-Altman's algebraic definition of ``fundamental forms" \cite[p.10]{Kleiman:Duality}, which presumably rooted in Grothendieck's work, we extend this to the higher order fundamental forms by defining the $k$-th Gauss map 
$$g_k:\bU_k\longrightarrow \nG(r,\mf{t}_{k})$$
on a suitable open subset $\bU_k$ of $X$ sending a point $x\in \bU_k$ to $\bT^{k}_x(X)$. The differential $dg_{k}$ and the fundamental form $\bF_{k+1}$ determine each other. Furthermore, locally given a vector field $\partial\in T_X$, one has a commutative diagram
$$\displaystyle \begin{CD} 
	P^{k}(L)^*   @>\alpha^*_{k}>> P^*_{k}\otimes L^*\\
	@VVV @VVdg_{k}(\partial)V \\
	S^{k}T_X\otimes L^*  @>\bF_{k+1,\partial}\otimes L^*>> R^*_{k}\otimes L^*
\end{CD}$$
where $*$ means taking the dual of a locally free sheaf (see Theorem \ref{thm:03} for details). In particular, at a generic point $x$, there is an induced linear map 
	$$f_{k+1,x}:T_x(X)\times \widehat{\bT}^{k}_x(X)\longrightarrow  V^*/\widehat{\bT}^{k}_x(X).$$
	This linear map gives us two interesting numerical invariants. Taking a generic line $\hat{y}\in \widehat{\bT}^{k}_x(X)$, the first invariant is the rank of the induced map $f_{k+1,x}(\hat{y}):T_x(X)\rightarrow V^*/\widehat{\bT}^{k}_x(X)$. It measures the defect of $k$-th order tangent variety (see Theorem \ref{thm:52} and  Remark \ref{rmk:41}). The second one, denoted by $\theta_k$, has less geometric intuition: it is defined to be the rank of the induced map $f_{k+1,x}(u): \widehat{\bT}^{k}_x(X)\rightarrow V^*/\widehat{\bT}^k_x(X)$ for a generic vector $u\in T_x(X)$, in symbols  
$$\theta_{k}=\rank f_{k+1,x}(u).$$
This number will be used in our vanishing theorem. 

We need one more numerical invariant to state our vanishing result. Inspired by the Terracini lemma, we define the number 
$$\delta_k=\dim \bT^k_x(X)\cap \bT^k_y(X), \text{ for generic points }x, y\in X.$$
Note that the number $\delta_1+1$ is the secant defect of $X$. 

\begin{theorem}\label{thm:vanishing} Let $X\subseteq\nP^r$ be a quasi-projective variety. For $k\geq 1$, if $\mf{t}_k=\theta_k+\delta_k$, then $\bF_{k+2,x}=0$ at a generic point $x\in X$.
\end{theorem}
The proof of the theorem is inspired by the work of Griffiths-Harris on third fundamental forms but using sheaves of principal parts developed here. We bypass the technical object called refined third fundamental form used in \cite{Griffiths:AGandLDG} to get the full generality. The theorem also gives a criterion when the variety is contained in a linear space since the vanishing of $\bF_{k+2}$ implies that $X$ is contained in a generic $(k+1)$-th projective tangent space.

We mention a special case of the above theorem. When $k=1$, the condition $\mf{t}_1=\theta_1+\delta_1$ is equivalent to the condition that the tangent variety $\tau(X)$ (the Zariski closure of the union of projective tangent space at nonsingular locus of $X$) equals the secant variety $\Sec(X)$. This forces $\tau(X)=\Tan(X)=\Sec(X)$ which is the degenerate case of Fulton-Hansen theorem. So the theorem claims that the third fundamental form $\bIII$ is zero, generalizing the results of Griffiths-Harris and Landsberg. We state it as a corollary below.

\begin{corollary}\label{thm:06}Let $X\subseteq\nP^r$ be a quasi-projective variety. If $\tau(X)=\Sec(X)$, then $\bIII_x=0$ at a generic point $x\in X$.
\end{corollary}
An immediate consequence is that $X$ in Corollary is contained in its generic second tangent space. We point out that the converse of the corollary is not true. For example, let $X$ be a generic projection of $3$-uple of $\nP^2$ from $\nP^9$ to $\nP^5$. Then $X$ has non-defective secant variety but $\bIII_x=0$ at generic points. The corollary gives a lower bound for the second fundamental form $\bII$ for a nondegenerate projective variety. 

\begin{corollary}\label{thm:07} Let $X\subseteq\nP^r$ be a nondegenerate projective variety with $\Tan(X)=\tau(X)$. Then at a generic point $x$, $\rank \bII_x\geq \min\{\codim X, \dim X\}$.
\end{corollary}

As a consequence, the generic second tangent space of the variety $X$ in Corollary \ref{thm:07} has dimension no less than $\min\{r,2\dim X\}$. This was claimed by F. Zak for nonsingular case in an online video of his seminar. However, we are unable to find a published reference. 

\bigskip
The paper is organized as follows. In section 2, we briefly review the notion of principal parts of line bundles and differential operators and then give the definition of fundamental forms. In section 3, we study higher order Gauss maps and projective geometry of higher order projective tangent spaces. In Section 4, we first study higher order tangent varieties and then give the proof of vanishing theorem of fundamental forms. 

\bigskip
\noindent{\bf Convention.} A variety is separated, reduced and irreducible of finite type over $\nC$.
Let $\sF$ be a coherent sheaf on a variety $X$ and $x\in X$ be a closed point with the residue field $k(x)$. We write $\sF_x(X)=\sF\otimes k(x)$. This notation is particularly applied to cotangent space $\Omega^1_x(X)$ or $T^*_x(X)$, Zariski tangent space $T_x(X)$, conormal space $N^*_x(X)$ and normal space $N_x(X)$. If $X$ is in a projective space $\nP^r$,  we write $\bT_x(X)$ to be the projective tangent space of $X$ at $x$ in $\nP^r$. 

\section{Fundamental forms of algebraic varieties}

\noindent In this section, we develop an algebraic approach to the theory of fundamental forms of algebraic varieties using principal parts of line bundles. As mentioned in Introduction, the modern approach was developed  by moving frames in differential geometry in the influential paper of Griffiths-Harris \cite{Griffiths:AGandLDG}.

Throughout the section, we assume $X\subseteq\nP^r$ be a quasi-projective variety of dimension $n\geq 1$, $L=\sO_X(1)$ and $V=H^0(\nP^r,\sO_{ \nP^{r}}(1))$.

\subsection{Principal parts of line bundles and differential operators}$ $
\medskip

We start by briefly reviewing the notion of principal parts of line bundles and differential operators. We refer the reader to \cite[Chapter 16]{Grothendieck:EGAIV1967} for details and full generality.

For an integer $k\geq 0$ the {\em sheaf of $k$-th order principal parts} of $L$ is defined to be  
$$P^k(L)=pr_{1,*}\Big((\sO_{X\times X}/I^{k+1}_{\Delta})\otimes pr^*_2 L\Big),$$
where $I_{\Delta}$ is the ideal sheaf of the diagonal $\Delta$ of $X\times X$ and $pr_1$ and $pr_2$ are projections of $X\times X$ to its components. We write $P^k_X$ for the sheaf $P^k(\sO_X)$ of $k$-th order principal parts of $\sO_X$. The  $\sO_X$-module structure of the sheaf $P^k_X(L)$ is induced by $pr_1$. The projection $pr_2$ induces the {\em  universal differential operator} (of order $\leq k$)
\begin{eqnarray*}
	d^k_L:L &\longrightarrow &P^k(L),
\end{eqnarray*} 
which is a  $\nC$-linear map.
In particular, for the $\sO_X$-algebra $P^k_X$, the universal differential operator $d^k:\sO_X\rightarrow P^k_X$ is a $\nC$-linear {\em algebra} map. The sheaf $P^k(L)$ is also a $P^k_X$-module and the differential operators satisfy the property that 
$$d^k_L(as)=d^k(a)d^k_L(s), \text{ for local sections }a\in \Gamma(U,\sO_X)\text{ and } s\in \Gamma(U,L).$$


Directly from definition there is a surjective {\em truncation map}
$\lambda_{k,k-1}:P^k(L)\rightarrow P^{k-1}(L)$ if $k\geq 1$, induced by the surjective morphism $\sO_{X\times X}/I^{k+1}_\Delta\rightarrow \sO_{X\times X}/I^{k}_\Delta$. The kernel sheaf of $\lambda_{k,k-1}$ is $pr_{1,*}(I^{k}_\Delta/I^{k+1}_{\Delta}\otimes pr_1^*L)$. If $X$ is {\em nonsingular}, $P^k_X(L)$ is locally free and the kernel sheaf of the truncation map $\lambda_{k,k-1}$ is $S^k\Omega^1_X\otimes L$. Thus in the case of $X$ nonsingular, one obtains a short exact sequence of locally free sheaves
\begin{eqnarray}\label{eq:08}
	0\longrightarrow S^k\Omega^1_X\otimes L\longrightarrow P^k(L)\longrightarrow P^{k-1}(L)\longrightarrow 0.
\end{eqnarray}
Applying $\sHom(\_,\sO_X)$ to the sequence above yields a short exact sequence 
\begin{equation}\label{eq:09}
	0\longrightarrow P^{k-1}(L)^*\longrightarrow P^{k}(L)^*\stackrel{pr_k}{\verylongrightarrow} S^kT_X\otimes L^*\longrightarrow 0,
\end{equation}
where $L^*=\sHom(L,\sO_X)$ and $P^k(L)^*=\sHom_{\sO_X}(P^k(L),\sO_X)$.

There is an evaluation morphism $e_L:V\otimes\sO_X\rightarrow L$ on global sections (we drop $\nC$ in $V\otimes_{\nC} \sO_X$ if there is no confusion arose). Pulling back  $e_L$ by $pr_2$ and then tensoring with $\sO_{X\times X}/I^{k+1}_\Delta$ gives a morphism $V\otimes \sO_{X\times X}\rightarrow (\sO_{X\times X}/I^{k+1}_\Delta)\otimes pr_2^*L$. Pushing down this map by $pr_1$ gives rise to a Taylor series map that we shall define below. 

\begin{definition}\label{def:05} For $k\geq 0$,  the ($k$-th)  {\em Taylor series map} 	$$\alpha_k:V\otimes \sO_X\longrightarrow P^k(L) $$
	is the $\sO_X$-homomorphism induced by the  evaluation morphism $V\otimes\sO_X\rightarrow L$.  Associated to the Taylor series map $\alpha_k$, define the sheaves
	$$R_k=(\Ker \alpha_k)\otimes L^*, \text{ and }P_k=(\uIm \alpha_k)\otimes L^*.$$
	Define the numbers $\mf{t}_k$ and $\mf{c}_k$ in the way that
	$$\mf{t}_k+1=\text{ the rank of the sheaf $P_k$},  \text{ and } \mf{c}_k=\text{the rank of the sheaf } R_k.$$  The truncation map $\lambda_{k,k-1}$ induces a surjective truncation map $\lambda'_{k,k-1}:P_k\otimes L\rightarrow P_{k-1}\otimes L$. Define the sheaf
	$$S_k=(\Ker \lambda'_{k,k-1})\otimes L^*\text{ if }k\geq 1.$$
\end{definition}

We frequently work at a generic point of $X$,  around which all  sheaves involved in Definition \ref{def:05} are locally free. To be more precise, we introduce the following open subsets $\bU_k$. Note that $\bU_0$ is the nonsingular locus of $X$ and $\bU_1=\bU_0$.

\begin{definition}\label{def:03} For $k\geq 0$, define $\bU_k$ to be the maximal open subset of $X$ contained in the nonsingular locus such that the quotient sheaf $P^i(L)/(P_i\otimes L)$ is locally free of constant rank on $\bU_k$ for all $i\leq k$.
\end{definition}

\begin{remark}\label{rmk:04} The Taylor series map $\alpha_k$ in Definition \ref{def:05} can be expressed locally on an open subset  $U$ of $X$ as
	$$\alpha_k(T\otimes a)=ad^k_L(\bar{t}), \text{ for }T\in V, \ a\in \Gamma(U,\sO_X).$$
	where $\bar{t}$ is the restriction of the section $T$ onto $U$ and $d^k_L$ is the universal differential operator.	Globally since $L$ is globally generated the Taylor series map $\alpha_0$ is surjective so that $P^0(L)=P_0\otimes L=L$ and $R_0\otimes L=\Omega^1_{\nP^r}|_X\otimes L$.  As a result, one gets the Euler sequence $$0\longrightarrow \Omega^1_{\nP^r}|_X\otimes L\longrightarrow V\otimes \sO_X\stackrel{\alpha_0}{\longrightarrow} L\longrightarrow 0.$$ If $X$ nonsingular, since $L$ is very ample, the Taylor series map $\alpha_1$ is surjective. Thus $P_1\otimes L=P^1(L)$, $R_1\otimes L=N_X^*\otimes L$, and one obtains the short exact sequence $$0\longrightarrow N^*_X\otimes L\longrightarrow V\otimes \sO_X\stackrel{\alpha_1}{\longrightarrow} P^1(L)\longrightarrow 0.$$
\end{remark}

%
%

We discuss local properties of Taylor series maps and principal parts of line bundles. It is convenient to replace $X$ by a suitable affine open subset of a nonsingular point. We shall apply differential operators locally and use multi-index notation such as ${\bf{p}}=(p_1,\ldots, p_n)$ following \cite[16.11]{Grothendieck:EGAIV1967}. So it is necessary to fix notation in the following

\begin{localsetup}\label{localsetup}	 Replacing $X$ by an {\em nonsingular affine open} set  such that there are local sections $s_1,\ldots,s_n\in\sO_X$ such that $\{ds_i\}$ is a basis for the free $\sO_X$-module $\Omega^1_{X}$. Define
	$$ds_i=d^ks_i-s_i\in P^k_{X},$$
	where  $d^k:\sO_{X}\rightarrow P^k_{X}$ is the  universal differential operator. For a multi-index ${\bf p}=(p_1,\ldots,p_n)$ define the product
	$$ds^{\bf p}=ds_1^{p_1}ds_2^{p_2}\ldots ds_n^{p_n}.$$
	The set $\{ds_i\}_{i=1}^n$ generates $\sO_X$-algebra $P^k_{X}$ and the set  $\{ds^{\bf p}\mid |{\bf p}|\leq k \}$  is a basis for the free $\sO_X$-module $P^k_{X}$. The dual module $D^k(\sO_X)=\sHom(P^k_X,\sO_X)$ has the dual basis $\{D_{\bf p}\mid |{\bf p}|\leq k\}$. In particular, the tangent sheaf $T_X$ has a basis $\{D_1,\ldots, D_n\}$.
\end{localsetup}

Let $x\in X$ be a nonsingular point. Work in Local Setup \ref{localsetup} on a suitable affine open neighborhood of $x$ so that $P^k(L)$ is a free $\sO_X$-module. Assume that 
\begin{equation}\label{tri}
	\text{there exists a section $t_0\in L$ trivializing $L$, i.e., $L=\sO_X\cdot t_0$}.
\end{equation} The set $\{ds^{\bf p}d^k_L(t_0)\mid |{\bf p}|\leq k\}$ is a basis for $P^k(L)$ as a module over $\sO_X$.  The dual module $P^k(L)^*=\sHom_{\sO_X}(P^k(L),\sO_X)$ has the dual basis $\{D_{\bf p}D^L_{t_0}\mid |{\bf p}|\leq k\}$, where $D^L_{t_0}:L\rightarrow \sO_X$ is a differential operator (of order zero) from $L$ to $\sO_X$ and $D_{\bf p} D^L_{t_0}$ is the composition of differential operators. One checks that  $D^L_{t_0}(at_0)=a$ and $D_{\bf p}D^L_{t_0}(a{t_0})=D_{\bf p}(a)$ for $a\in \sO_X$.

%

%
For an element $e\in \sO_X$, the universal differential operator $d^k:\sO_X\rightarrow P^k_X$ has the expression $d^k(e)=\sum_{0\leq |{\bf p}|\leq k}D_{\bf p}(e)ds^{\bf p}$.
Since $L$ is trivialized by the section $t_0$, for any section $T\in V$, the restriction of $T$ onto $L$ can be uniquely written as $t\cdot t_0$ for an element $t\in \sO_X$. In this case,  for an element $T\otimes a\in V\otimes\sO_X$ with $a\in \sO_X$, 
\begin{equation}\label{localexpress}
	\alpha_k(T\otimes a)=ad^k_L(t\cdot t_0)=ad^k(t)d^k_L(t_0)=a\sum_{0\leq |{\bf p}|\leq k}D_{\bf p}(t)ds^{\bf p}d^k_L(t_0).
\end{equation}
We point out that one can obtain an isomorphism $P^k(L)\cong P^k_X$ by sending $d^k_L(t_0)$ to $d^k(1)$. Under this isomorphism the Taylor series map can be considered as a map $\alpha_k:V\otimes \sO_X\rightarrow P^k_X$ and one can drop $d^k_L({t}_0)$ in the above local expression (\ref{localexpress}). However, the trivialization (\ref{tri}) of $L$ is not canonical so sometimes we prefer to bring $d^k_L$ in the calculation. 

For a vector field $\partial\in T_X$, there is an induced $\nC$-homomorphism  (we use the same notation)
$$\partial:P^{k-1}(L)^*\longrightarrow P^k(L)^*$$
sending a differential operator $D\in P^{k-1}(L)^*$ to the composition $\partial D:=\partial\circ D\in P^k(L)^*$ (\cite[Proposition 16.8.9]{Grothendieck:EGAIV1967}). By abuse of notation, we may call it the  derivative of $P^{k-1}(L)^*$ into $P^k(L)^*$. In terms of the dual basis $\{D_{\bf p}D^L_{t_0}\mid |{\bf p}|\leq k-1\}$ of $P^{k-1}(L)^*$, the map $\partial$  can be written as
$$\partial(aD_{\bf p}D^L_{t_0})= \partial (a)D_{\bf p}D^L_{t_0}+a\partial D_{\bf p}D^L_{t_0}, \text{ for }a\in \sO_X.$$
Since the projection map $	pr_k:P^k(L)^*\rightarrow S^kT_X\otimes L^*$ in  (\ref{eq:09}) is defined by 
\begin{eqnarray*}
	\sum_{|{\bf p}|\leq k}a_{\bf p}D_{\bf p}D^L_{t_0}&\mapsto & \sum_{|{\bf p}|=k}a_{\bf p}D_{\bf p}D^L_{t_0}, \text{ where }a_{\bf p}\in \sO_X,
\end{eqnarray*}
one immediately checks that the vector field $\partial$ induces an $\sO_X$-homomorphism 
$$\partial:	S^{k-1}T_X\otimes L^* \longrightarrow S^kT_X\otimes L^*,$$
sending $D_{\bf p}D^L_{t_0}$ to $\partial D_{\bf p}D^L_{t_0}$. Finally, as $\partial:\sO_X\rightarrow \sO_X$ is a $\nC$-homomorphism, tensoring with $V^*$ over $\nC$ induces an $\nC$-homomorphism  
$$\partial:V^*\otimes \sO_X\longrightarrow  V^*\otimes \sO_X$$
sending an element $T\otimes a$ to $T\otimes \partial(a)$. We also call it the derivative on $V^*\otimes \sO_X$. Taking a section $B\in V^*\otimes \sO_X$, the evaluation of $B$ at the point $x$ is a vector $B(x)\in V^*$. Applying the derivative $\partial$ to $B$ yields an element  $\partial B$ in $V^*\otimes \sO_X$ and thus a vector $(\partial B)(x)\in V^*$.

\begin{proposition}\label{p:05} Work in  Local Setup \ref{localsetup}. 
	\begin{enumerate}
		\item For a vector field $\partial \in T_X$, we have a commutative diagram 
		$$\begin{CD}
			S^{k-1}T_X\otimes L^* @<pr_{k-1}<<	P^{k-1}(L)^* @>\alpha^*_{k-1}>> V^*\otimes \sO_X \\
			@V\partial VV	@V\partial VV @V\partial VV\\
			S^kT_X\otimes L^* @<pr_k<<	P^k(L)^* @>\alpha^*_{k}>> V^*\otimes \sO_X
		\end{CD}.$$
		\item The $\sO_X$-module $P^k(L)^*$ is generated by the submodule $P^{k-1}(L)^*$ and its derivatives of all $\partial \in T_X$, i.e., 
		$$P^k(L)^*=P^{k-1}(L)^*+\sum_{\partial \in T_X} \partial P^{k-1}(L)^*.$$
		\item As an $\sO_X$-submodule of $V^*\otimes \sO_X$, the image of $\alpha_k^*$ is generated by the image of $\alpha_{k-1}^*$ and its derivatives of all $\partial \in T_X$, i.e., 
		$$ \uIm(\alpha^*_k)=\uIm(\alpha^*_{k-1})+\sum_{\partial\in T_X}\partial (\uIm(\alpha^*_{k-1})).$$
		
	\end{enumerate}
\end{proposition}

\begin{proof} (1) The left-hand-side commutative square can be easily checked by definition. So we only need to prove the right-hand-side commutative square. Let $\{T_i\}$ be a basis of the vector space $V$ and let $\{T^*_i\}$ be the dual basis for $V^*$ so that $V^*\otimes \sO_X=\sO_XT^*_0\oplus \sO_XT^*_1\oplus \ldots\oplus \sO_XT^*_r$. Recall that $L$ is trivialized  as $L=\sO_X\cdot t_0$ by a local section $t_0$, the restriction of $T_i$ to $L$ can be uniquely written as $\bar{t}_i\cdot t_0$ for an element $\bar{t}_i\in \sO_X$.  The dual of the Taylor series map is given by 
	\begin{eqnarray*}
		\alpha_k^*:P^{k}(L)^*&\longrightarrow &V^*\otimes \sO_X=\sO_XT^*_0\oplus \sO_XT^*_1\oplus \ldots\oplus \sO_XT^*_r\\
		D_{\bf p}D^L_{t_0} &\mapsto & (D_{\bf p}(\bar{t}_0), \ldots, D_{\bf p}(\bar{t}_r))
	\end{eqnarray*}
	where  $r$-tuple $(D_{\bf p}(\bar{t}_0), \ldots, D_{\bf p}(\bar{t}_r))$ represents the element $D_{\bf p}(\bar{t}_0)T^*_0+\ldots+D_{\bf p}(\bar{t}_r)T^*_r$ in $V^*\otimes \sO_X$. Now the the desired commutativity of the diagram can be checked directly by definition. 
	
	(2) and (3). Keep using the notation in (1). For each multi-index $\bf p$, we write $e_{\bf p}=(D_{\bf p}(\bar{t}_0), \ldots, D_{\bf p}(\bar{t}_r))$ as an element of $V^*\otimes \sO_X$. 
	We see $\uIm(\alpha^*_k)$ is generated by $\{ e_{\bf p}\mid |{\bf p}|\leq k\}$.	Similarly, 
	$\uIm(\alpha^*_{k-1})$ is generated by $\{ e_{\bf p}\mid |{\bf p}|\leq k-1\}$.
	Note that $T_X=\sO_XD_1\oplus \ldots \oplus \sO_XD_n$ and $D_i(e_{\bf p})=(D_i\circ D_{\bf p}(\bar{t}_0),\ldots, D_i\circ D_{\bf p}(\bar{t}_r))=e_{\bf p+i}$. Thus the results can be easily checked by definition.  
\end{proof}

\begin{proposition}\label{p:20} The rank sequence $\{{\mf c}_1,\mf{c}_2,\cdots,\}$ of the sheaves $R_i$ is strictly decreasing to a stable number, i.e., $\mf{c}_1>\mf{c}_2>\ldots>\mf{c}_m=\mf{c}_{m+1}=\ldots$ such that \begin{enumerate}
		\item either $\mf{c}_m=0$, $R_m=0$, and $X$ is not contained in any sub linear space of $\nP^r$, or
		\item there exists a subspace $W\subseteq V$ of dimension $\mf{c}_m$ such that $R_m\otimes L=W\otimes \sO_X$ contained in all $R_i\otimes L$ for $i\geq 0$, and $X$ is contained in the linear space $\nP(V/W)$.
	\end{enumerate}
\end{proposition}
\begin{proof} We show first that the sequence $\{\mf{t}_1,\mf{t}_2,\ldots\}$ associated to the sheaves $P_i$ is strictly increasing to a stable number. Observe first that it always increases and has a upper bound $r+1$. So it suffices to show that if for some $m\geq 1$, $\mf{t}_m=\mf{t}_{m+1}$ then $\mf{t}_i=\mf{t}_m$ for all $i\geq m$. To see this, work on the open subset $\bU_{m+1}$. The sheaves $P_{m}$ and $P_{m+1}$ are all locally free of rank $\mf{t}_m+1$. So shrinking $\bU_{m+1}$ if necessary, we assume $P_{m}=P_{m+1}$. By Proposition \ref{p:05}, for any vector field $\partial\in T_X$, and any section $B\in P_{m}$, we have $\partial(B)\in P_{m}$. Thus by Proposition \ref{p:05} again, we see on the open subset $\bU_{m+i}$ (or its suitable open subset) $P_m=P_{m+i}$ for all $i\geq 0$. This proves the sequence strictly increases to a stable number. As a consequence, the rank sequence $\{\mf{c}_1,\mf{c}_2,\ldots, \}$ is strictly decreasing to a stable number $\mf{c}_m$ as claimed. If $\mf{c}_m=0$, then since the sheaf $R_m$ is a torsion-free sheaf, it must be zero. In the sequel, we assume $\mf{c}_m\neq 0$.  As we have countably many open subsets $\bU_i$, $i\geq 0$, there exists a closed point $x$ in every open subset $\bU_i$. Thus we have $R_m\otimes k(x)=R_{m+1}\otimes k(x)=\ldots$ and there is a short exact sequence $0\rightarrow R_m\otimes k(x)\rightarrow V\rightarrow P_m\otimes k(x)\rightarrow 0$.
	Let $W=R_m\otimes k(x)$ as a subspace of $V$. Thus for a section $T\in W$, we have $\alpha_{k,x}(T)=0$ for all $k\geq 0$. As $X$ is irreducible, so $X$ is contained in the linear space defined by $W$. Thus globally, the trivial sheaf $W\otimes \sO_X$ is contained in every sheaf $R_i\otimes L$ for $i\geq 0$. Furthermore, we have the following diagram 
	$$\displaystyle \begin{CD} 
		0 @>>> R_i\otimes L  @>>> V\otimes \sO_X @>\alpha_m>> P^m(L) @.\\
		@. @VVV @VVV @|\\
		0 @>>> R'_m\otimes L  @>>> (V/W)\otimes \sO_X @>\alpha'_m>>P^m(L) @.
	\end{CD}$$
	The Taylor series maps $\alpha_m$ and $\alpha'_m$ have the same image. By the Snake lemma, we obtain a short exact sequence $0\rightarrow W\otimes \sO_X\rightarrow R_m\otimes L\rightarrow R'_m\otimes L\rightarrow 0$.
	Note that $W\otimes \sO_X$ and $R_m\otimes L$ have the same rank and the sheaf $R'_m\otimes L$ is torsion free. Hence we conclude that $R'_m\otimes L=0$ and therefore $W\otimes \sO_X=R_m\otimes L$.
\end{proof}

\subsection{Definition of fundamental forms}$ $

\medskip
The fundamental form that we will define is induced by the Taylor series map.  Since we exclusively work on the nonsingular locus of $X$, without loss of generality, we may assume $X$ is nonsingular. Restrict $\alpha_k$ to the subsheaf $R_{k-1}\otimes L$ and one checks by Snake lemma that the image is landing in the kernel sheaf of the truncation map $\lambda_{k,k-1}$, which is $S^{k}\Omega^1_X\otimes L$. Thus we obtain an induced map
$$\alpha_k|_{R_{k-1}\otimes L}:R_{k-1}\otimes L\longrightarrow S^k\Omega^1_X\otimes L.$$
\begin{definition}\label{def:02} Assume $X$ is nonsingular. For $k\geq 1$, the $k$-th {\em twisted fundamental form} is defined to be the map $\bF^*_k\otimes L:=\alpha_k|_{R_{k-1}\otimes L}$ and the $k$-th {\em fundamental form} is the induced map $$\bF_k^*:R_{k-1}\longrightarrow S^k\Omega^1_X$$
	obtained  by twisting $\bF^*_k\otimes L$ by $L^*$. We also write $\bII^*$ and $\bIII^*$ for the {\em second and third fundamental form} respectively. We denote by $\bF_k$ the dual map obtained by applying $\sHom_{\sO_X}(\_,\sO_X)$ to $\bF^*_k$ and similarly for $\bII$ and $\bIII$. The morphisms $\bF_k, \bII, \bIII$, etc, are also called the fundamental forms and $\bF_k\otimes L^*$ are also called the twisted fundamental forms.
\end{definition}

\begin{definition}\label{def:06} Let $x\in X$ be a nonsingular point. Tensoring with the residual field $k(x)$ to the $k$-th fundamental form $\bF^*_k$ yields a morphism $$\bF^*_{k,x}: R_{k-1,x}(X)\longrightarrow S^k\Omega^1_x(X)$$
	on vector spaces. Similar notations work for $\bF_{k,x}, \bII^*_x, \bII_x$, and twisted $\bF_{k,x}\otimes L^*$ etc. 
	
\end{definition}
\begin{remark} The notations in Definition \ref{def:02} and \ref{def:06} are chosen to be consistent with the ones appeared in the early work such as \cite{Griffiths:AGandLDG} and \cite{Landsberg:DegSecTan}. Using the sheaves $P_k$ and $S_k$ in Definition \ref{def:05}, we can form  the following diagram 
	\begin{equation}\label{eq:07}
		\displaystyle \begin{CD} 
			@.  @. @.0\\
			@. @. @. @VVV\\
			@.@. @. S_{k}\otimes L  @.\\
			@. @. @. @VVV\\
			0 @>>> R_{k}\otimes L @>>> V\otimes \sO_X @>\alpha_k>>P_{k}\otimes L @>>>0\\
			@. @VVV @| @VVV\\
			0@>>> R_{k-1}\otimes L@>>> V\otimes \sO_X @>>> P_{k-1}\otimes L@>>> 0\\
			@. @. @. @VVV\\
			@.  @.  @.0. @.\\
		\end{CD}
	\end{equation}
We see that the sheaf $S_k\otimes L$ is the image sheaf of the twisted fundamental form $\bF_k^*\otimes L$  and the sheaf $S_k$ is the image sheaf of the fundamental form $\bF_k^*$. 
The diagram works for singular case so one can define the (twisted) fundamental forms on arbitrary quasi-projective variety. 
\end{remark}

Directly from definition, we can give a local formula for fundamental forms  around a nonsingular point. Recall that if $X$ is nonsingular, then the sheaf $R_2$ is the conormal sheaf $N^*_X$ and all $R_k$ with $k\geq 2$ are subsheaves of $N^*_X$. Let $I_X$ be the defining ideal sheaf of $X$, then by definition the conormal sheaf $N^*_X=I_X/I^2_X$. So an element of $N^*_X$ can be written as a quotient class $\bar{f}$ for an element $f\in I_X$. 
\begin{proposition}\label{p:03} Let $x\in X$ be a nonsingular point and work in Local Setup \ref{localsetup} on an affine open neighborhood of $x$. Let $\{T_0,T_1,\ldots, T_r\}$ be a basis of the vector space $V$ and let $U\subseteq \nP^r$ be an suitable affine open subset containing $X$ such that  
	\begin{enumerate}
		\item the restriction $t_0=T_0|_U$ gives a local trivialization $\sO_U(1)|_U=\sO_U\cdot t_0$;
		\item for $i=1,\ldots,r$, $T_i|_U=t_it_0$ with $t_i\in \sO_U$ such that $\{dt_1,\ldots, dt_r\}$ is a basis  for $\Omega^1_{\nP^r}|_U$. 
	\end{enumerate}
	Write $\bar{t}_i$ the restriction of $t_i$ in the ring $\sO_X$ and $\bar{t}_0$ the restriction of $t_0$ in $L$. For $k\geq 2$ and an element $\bar{f}\otimes \bar{t}_0\in N^*_X\otimes L=I_X/I_X^2\otimes L$, where $f\in I_X$, the Taylor series map $\alpha_k:N^*_X\otimes L\rightarrow P^k(L)$ is given by 
	$$\alpha_k(\bar{f}\otimes \bar{t}_0)=\Big (\sum_{i=1}^r \frac{\overline{\ud f}}{\ud t_i}(\sum_{2\leq |{\bf p}|\leq k}D_{\bf p}(\bar{t}_i)ds^{\bf p})\Big )d^k_L(\bar{t}_0),$$
	where $\frac{\overline{\ud f}}{\ud t_i}$ is the restriction of $\frac{\ud f}{\ud t_i}$ to the ring $\sO_X$.
	
\end{proposition}

\begin{proof} The Taylor series map $\alpha_{k}:V\otimes \sO_{X}\rightarrow P^k(L)$	is defined by $\alpha_k(T\otimes a)=ad^k_L(\bar{t})$ for $T\in V$, $a\in \sO_X$,	where $\bar{t}$ is the restriction of $T$ to $L$.
	On the other hand, the inclusion map $\rho: N^*_X\otimes L\rightarrow V\otimes \sO_X$ is sending $\bar{f}\otimes \bar{t}_0$ to $\sum_{i=1}^r \frac{\overline{\ud f}}{\ud t_i}(T_i\otimes 1-T_0\otimes \bar{t}_i)$.
	Thus we calculate directly that $\alpha_k(\rho(\bar{f}\otimes \bar{t}_0))=\Big(\sum_{i=1}^r \frac{\overline{\ud f}}{\ud t_i}(d^k\bar{t}_i-\bar{t}_i)\Big)d^k_L(\bar{t}_0)$, which equals $\Big(\sum_{i=1}^r \frac{\overline{\ud f}}{\ud t_i}(\sum_{1\leq |{\bf p}|\leq k}D_{\bf p}(\bar{t}_i)ds^{\bf p})\Big)d^k_L(\bar{t}_0)$ since $d^k\bar{t}_i=\sum_{0\leq |{\bf p}|\leq m}D_{\bf p}(\bar{t}_i)ds^{\bf p}$ and $D_{\bf 0}(\bar{t}_i)=\bar{t}_i$.	But observe that if $|{\bf p}|=1$, the  term $\sum_{i=1}^r \frac{\overline{\ud f}}{\ud t_i}(\sum_{|{\bf p}|=1}D_{\bf p}(\bar{t}_i)ds^{\bf p}) =d\bar{f}=0$. Hence we obtain the desired formula in the proposition. 
\end{proof}

\begin{corollary}\label{p:04} With notation and assumption in Proposition \ref{p:03}, for an element $\bar{f}\in N^*_X=I_X/I_X^2$ represented by $f\in I_X$, one has 
	\begin{enumerate}
		\item the following are equivalent 
		$$\bar{f}\in R_{k-1}\ \Longleftrightarrow\ \alpha_{k-1}(\bar{f}\otimes \bar{t}_0)=0 \ \Longleftrightarrow\ \bF^*_j(\bar{f})=0 \text{ for }1\leq j\leq k-1.$$
		\item 	The $k$-th fundamental form $\bF^*_k:R_{k-1}\rightarrow S^k\Omega^1_X$ has the local form 
		$$\bF^*_k(\bar{f})=\sum_{i=1}^r \frac{\overline{\ud f}}{\ud t_i}(\sum_{|{\bf p}|=k}D_{\bf p}(\bar{t}_i)ds^{\bf p}), \text{ for } \bar{f}\in R_{k-1},$$
		which does not depend on the choice of the trivialization of $L$ described in Proposition \ref{p:03}.
		\item 	Tensoring with the residue field $k(x)$, the fundamental form $\bF^*_{k,x}:R_{k-1,x}(X)\rightarrow S^k\Omega^1_x(X)$
		has the expression 
		$$\bF^*_{k,x}(\bar{f}\otimes 1)=\sum_{i=1}^r \frac{\overline{\ud f}}{\ud t_i}(x)(\sum_{|{\bf p}|=k}D_{{\bf p},x}(\bar{t}_i)S^{\bf p}),$$
		where  $\frac{\overline{\ud f}}{\ud t_i}(x)$ and $D_{{\bf p},x}(\bar{t}_i)$ are the evaluation of $\frac{\overline{\ud f}}{\ud t_i}$ and $D_{{\bf p}}(\bar{t}_i)$ at $x$ respectively and $S^{\bf p}=S_1^{p_1}\cdot\ldots\cdot S_n^{p_n}$ with $S_i=ds_i\otimes 1$.
	\end{enumerate}

\end{corollary}\label{p:19}
\begin{proof} The equivalence statements in (1) are directly from the definition of fundamental forms. To prove (2), by Proposition \ref{p:03}, for $\bar{f}\in R_{k-1}$, we have 
	$\alpha_k(\bar{f}\otimes \bar{t}_0)=(\sum_{i=1}^r \frac{\overline{\ud f}}{\ud t_i}(\sum_{|{\bf p}|=k}D_{\bf p}(\bar{t}_i)ds^{\bf p}))d^k_L(\bar{t}_0)$
	where recall $L=\sO_X\cdot \bar{t}_0$. The image lands in $S^k\Omega_X^1\otimes L$ and the inclusion $S^k\Omega_X^1\otimes L\rightarrow P^k(L)$ is given by the map sending $ds^{\bf p}\otimes \bar{t}_0\mapsto ds^{\bf p}d^k_L(\bar{t}_0)$. Thus 
	$$\alpha_k(\bar{f}\otimes \bar{t}_0)=(\sum_{i=1}^r \frac{\overline{\ud f}}{\ud t_i}(\sum_{|{\bf p}|=k}D_{\bf p}(\bar{t}_i)ds^{\bf p}))\otimes \bar{t}_0\in S^k\Omega_X^1\otimes L.$$
	Hence tensoring with $L^*$, we get the desired formula in (2). The formula is independent on the trivialization of $L$. Indeed if $L=\sO_X\bar{t}_0'$ is another trivialization so that one has the corresponding $t'_i\in \sO_X$ involved in the formula. Then each $t_i=u_it'_i$ for some unit $u\in \sO_X$. Now use the fact $\bF^*_j(\bar{f})=0 \text{ for }1\leq j\leq k-1$ from (1) and Leibniz formula $D_{\bf p}(t)=D_{\bf p}(ut')=\sum_{{\bf q}\leq {\bf p}}{{\bf p}\choose {\bf q}}D_{\bf q}(u)D_{{\bf p}-{\bf q}}(t')$ to check immediately. Statement (3) is a directly consequence of (2).
\end{proof}


\begin{remark}We denote by $|\bF_{k,x}|$ the projectivization of the image of $\bF^*_{k,x}$ and consider it as a linear system on the projective space $\nP(T_x(X))$, where $T_x(X)$ is the Zariski tangent space of $X$ at $x$. In this way  $S_i=ds_i\otimes 1$ are considered as variables so that an element in $|\bF_{k,x}|$ is a degree $k$ homogeneous polynomial in $S_i$.
\end{remark}

\begin{remark} Working on the open subset $\bU_{k-1}$, the sheaf $R_{k-1}$ is locally free. Following from Corollary \ref{p:19} (3) one obtains an expression for the fundamental form $\bF_{k,x}:S^kT_x(X)\rightarrow R^*_{k-1,x}(X)$ and it shows that $\bF_{k,x}$ is a symmetric multi-linear map on $T_x(X)$. 
	
\end{remark}

\begin{remark} Using the local formula, one can  give an algebraic proof for a fundamental result due to Cartan  known in differential geometry (see \cite[4.2]{Landsberg:AlgGeoProjDiffGeo} for details) that the Jacobian system of $|\bF_{k,x}|$ is contained in the system $|\bF_{k-1,x}|$ at a generic point $x$. The local formula can also be used to show a well-known result that if $\bF_{k,x}=0$ at generic point $x$ then $X$ is contained in $\bT^{k-1}_x(X)$ (Definition \ref{def:39}).
\end{remark}

\section{Higher order Gauss map and projective geometry}

\noindent In this section, we discuss the geometry of fundamental forms by relating them to the differentials of higher order Gauss maps. Throughout the section, recall that $X\subseteq\nP^r$ is a quasi-projective variety of dimension $n\geq 1$, $L=\sO_X(1)$ and $V=H^0(\nP^r,\sO_{ \nP^{r}}(1))$.
\subsection{Higher order Gauss maps}$ $
\medskip

Consider the twisted fundamental form $\bF_k\otimes L^*$ on the open subset $\bU_{k-1}$. We can construct the following commutative diagram
\begin{equation}\label{eq:30}
	\begin{CD}
	\xymatrix{
		P^{k-1}(L)^* \ar@{->>}[r]^{\alpha^*_{k-1}} \ar@{^{(}->}[d] & P^*_{k-1}\otimes L^*   \ar@{^{(}->}[d]\ar@{^{(}->}[r] & V^*\otimes \sO_X \ar@{=}[d]  \ar@{->>}[r] & R^*_{k-1}\otimes L^*   \\
		P^{k}(L)^*\ar[r]^{\alpha^*_k} \ar@{->>}[d]& P^*_k\otimes L^*  \ar@{^{(}->}[r] \ar@{->>}[d]& V^*\otimes \sO_X  \\
		S^kT_X\otimes L^* \ar[r]^{\bF_k\otimes L^*}	& S^*_k\otimes L^*
	}\end{CD}
\end{equation}
Since $P_{k-1}$ is locally free on $\bU_{k-1}$, both vertical sequences on the left-hand are exact (the sheaves $P^*_k$ and $S^*_k$ are not necessarily locally free on $\bU_{k-1}$ but locally free on $\bU_k$). Thus $S^*_k\otimes L$ is a subsheaf of $R^*_{k-1}\otimes L^*$ containing the image of $\bF_k\otimes L^*$. Furthermore, if we work on the open subset $\bU_k$, then the map $\alpha^*_k$ becomes surjective and as a consequence  $S^*_k\otimes L^*$ equals the image sheaf of $\bF_k\otimes L^*$.

Work in Local Setup \ref{localsetup} on an affine open subset contained  in $\bU_{k-1}$. Recall that a vector field $\partial \in T_X$ induces $\nC$-homomorphisms  $\partial: P^{k-1}(L)^*\rightarrow P^k(L)^*$ and $\partial:V^*\otimes \sO_X\rightarrow V^*\otimes \sO_X$. By Proposition \ref{p:05} (3), we obtain an induced $\nC$-homomorphism $\partial:P^*_{k-1}\otimes L^*\rightarrow P^*_{k}\otimes L^*$ and  pass it to quotients to induce an $\sO_X$-homomorphism  $\partial:S^*_{k-1}\otimes L^*\rightarrow S^*_k\otimes L^*$, 
which fits into the following commutative diagram 
\begin{equation}\label{eq:81}
	\begin{CD}
		S^{k-1}T_X\otimes L^*@>\bF_{k-1}\otimes L^*>> S^*_{k-1}\otimes L^* @.\subseteq R^*_{k-2}\otimes L^*\\
		@V\partial VV @V\partial VV\\
		S^kT_X\otimes L^* @>\bF_k\otimes L^*>> S^*_k\otimes L^*@.\subseteq R^*_{k-1}\otimes L^*.
	\end{CD}
\end{equation}
We define 
\begin{equation}\label{eq:10}
	\bF_{k,\partial}\otimes L^*: S^{k-1}T_X\otimes L^*\longrightarrow R^*_{k-1}\otimes L^*
\end{equation}
to be the composition $(\bF_k\otimes L^*)\circ \partial $ in the diagram, which is a $\sO_X$-homomorphism.

\begin{definition}[{\bf $k$-th tangent space}]\label{def:39} The $k$-th projective tangent space $\bT^k_x(X)$ of $X$ at a point $x$ is defined to be $\bT^k_x(X)=\nP(P_k\otimes L\otimes k(x))$ which is a linear space in $\nP^r$. We write $\widehat{\bT}_x^{k}(X)$ to be the deprojectivization of $\bT_x(X)$ in $V^*$.
\end{definition}

\begin{definition}[{\bf $k$-th Gauss map}] Let $\mf{t}_k$ be the dimension of the $k$-th tangent space $\bT_x^k(X)$ for a point  $x\in\bU_k$. Define the {\em $k$-th Gauss map} $g_k:\bU_k\longrightarrow \nG(\nP^r, \mf{t}_k)$	by sending $x$ to $\bT^k_x(X)$.
\end{definition}

For $k\geq 2$,  the $(k-1)$-th Gauss map $g_{k-1}$ on the open set $\bU_{k-1}$ is determined by the short exact sequence
\begin{equation}\label{eq:03}
	0\longrightarrow R_{k-1}\otimes L\longrightarrow V\otimes \sO_X\longrightarrow P_{k-1}\otimes L\longrightarrow 0.
\end{equation} So the pullback of the tangent sheaf of the Grassmannian is $g_{k-1}^*T_{\nG}=\sHom(P^*_{k-1}\otimes L^*,R^*_{k-1}\otimes L^*)$. Locally the differential $dg_{k-1}$ sends a vector field $\partial\in T_X$ to a morphism $dg_{k-1}(\partial):P^*_{k-1}\otimes L^*\rightarrow R^*_{k-1}\otimes L^*$.

\begin{theorem}\label{thm:03} For $k\geq 2$, on the open set $\bU_{k-1}$ of $X$, consider the differential map 
	$$dg_{k-1}:T_X\longrightarrow g_{k-1}^*T_{\nG}=\sHom(P_{k-1}^*\otimes L^*, R^*_{k-1}\otimes L^*)$$
	of the $(k-1)$-th Gauss map.
	\begin{enumerate}
		\item The $k$-th fundamental form $\bF^*_{k}:R_{k-1}\rightarrow S^{k}\Omega^1_X$ is determined by the differential $dg_{k-1}$ and vice versa.
		\item With the Local Setup \ref{localsetup} on an affine open subset contained in $\bU_{k-1}$, for a vector field $\partial \in T_X$, the induced map 
		$$\bF_{k,\partial}\otimes L^*:S^{k-1}T_X\otimes L^*\longrightarrow R^*_{k-1}\otimes L^*$$
		defined in (\ref{eq:10}) fits into a commutative diagram 
		$$\displaystyle \begin{CD} 
			P^{k-1}(L)^*   @>\alpha^*_{k-1}>> P^*_{k-1}\otimes L^*\\
			@Vpr_{k-1}VV @VVdg_{k-1}(\partial)V \\
			S^{k-1}T_X\otimes L^*  @>\bF_{k,\partial}\otimes L^*>> R^*_{k-1}\otimes L^*.\\
		\end{CD}$$
		
	\end{enumerate}

\end{theorem}
\begin{proof} (1) The question is local so we work in Local Setup \ref{localsetup} on an affine open subset contained in $\bU_{k-1}$. Equivalently, we show the statement for the twisted fundamental form $\bF^*_k\otimes L$. Recall that $\bF^*_k\otimes L=\alpha_k|_{R_{k-1}\otimes L}:R_{k-1}\otimes L\rightarrow P^k(L)$, the restriction of the Taylor series map $\alpha_k$ onto the sub-module $R_{k-1}\otimes L$. 	
	Since on $\bU_{k-1}$, both $R_{k-1}$ and $P^k(L)$ are locally free,  $\bF^*_k\otimes L$ is determined by its dual map $\bF_k\otimes L^*$. As indicated in the diagram 
	$$\displaystyle \begin{CD} 
		0 @>>> P^*_{k-1}\otimes L^*  @>>> V^*\otimes \sO_X @>\pi_{k-1}>>R^*_{k-1}\otimes L^* @>>> 0\\
		@. @. @| @.\\
		@. P^{k}(L)^*  @>\alpha^*_k>> V^*\otimes \sO_X, @. @.\\
	\end{CD}$$
	$\bF_k\otimes L^*=\pi_{k-1}\circ \alpha^*_k$. Recall that the differential $dg_{k-1}:T_X\rightarrow \sHom(P_{k-1}^*\otimes L^*, R^*_{k-1}\otimes L^*)$ is defined by sending a vector field $\partial\in T_X$ to a map $dg_{k-1}(\partial):P^*_{k-1}\otimes L^*\rightarrow R^*_{k-1}\otimes L^*$  which by definition  $dg_{k-1}(\partial)(B)=\pi_{k-1}(\partial(B))$ for a section $B\in P^*_{k-1}\otimes L$. But Proposition \ref{p:05} tells us that $\partial (B)\in \uIm(\alpha^*_k)$. Hence there exists a section $A\in P^{k}(L)^*$   such that $\alpha^*_k(A)=\partial(B)$. This implies that $dg_{k-1}(\partial)(B)=\bF_k\otimes L^*(A)$, which means that $dg_{k-1}$ is determined by $\bF_k\otimes L^*$. Conversely, Proposition \ref{p:05} says that $\uIm(\alpha^*_k)$ is generated by $P^*_{k-1}\otimes L^*$ and $\partial (B)=dg_{k-1}(\partial)(B)$ for $B\in P^*_{k-1}\otimes L^*$ and $\partial\in T_X$. Thus $\bF_k\otimes L^*$ is determined by $dg_{k-1}$.
	
	(2) 
	The composition $dg_{k-1}(\partial)\circ \alpha^*_{k-1}$ maps an section $B\in P^{k-1}(L)^*$ to $\pi_{k-1}(\partial(\alpha^*_{k-1}(B)))$. By Proposition \ref{p:05}, $dg_{k-1}(\partial)\circ \alpha^*_{k-1}$ maps $P^{k-2}(L)^*$ to zero and thus factor through the the quotient $S^{k-1}T_X\otimes L^*$ to have an induced map $\psi$ as indicated in the following diagram 
	$$\xymatrix{
		P^{k-1}(L)^*\ar[d]_{pr_{k-1}} \ar[r]^-{\partial\circ \alpha^*_{k-1}}&	V^*\otimes \sO_X \ar[r]^{\pi_{k-1}} & R^*_{k-1}\otimes L^* \\
		S^{k-1}T_X \otimes L^* \ar[urr]_{\psi} &	 &  
	}$$
	So it suffices to check $\psi=\bF_{k,\partial}\otimes L^*$. To see this, we can form the following diagram
	$$\xymatrix{
		& P^{k-1}(L)^*\ar[r]^{\alpha^*_{k-1}}\ar[d]^{\partial} \ar[dl]_{pr_{k-1}} &	P^*_{k-1}\otimes L^* \ar[d]^{\partial}&  \\
		S^{k-1}T_X\otimes L^* \ar[dr]^{\partial}&		P^k(L)^* \ar[r]\ar[d]^{pr_k}\ar[r]^{\alpha^*_k}& P^*_k\otimes L^*	\ar[d]^{\pi:=\pi_{k-1}|_{P^*_k\otimes L^*}} \ar@{^{(}->}[r] & V^*\otimes \sO_X \ar[r]^{\pi_{k-1}}& R^*_{k-1}\otimes L^* \\
		& S^kT\otimes L^*\ar[r]^-{\bF_k\otimes L^*} & S^*_k\otimes L^*
	}$$
	The composition $\pi_{k-1}\circ \partial\circ \alpha^*_{k-1}=\pi\circ \partial \circ \alpha^*_{k-1}$. By commutativity, we see $\partial\circ \alpha^*_{k-1}=\alpha^*_k\circ \partial$ (Proposition \ref{p:05}), $\pi\circ \alpha^*_k=(\bF_k\otimes L^*)\circ pr_{k}$ (diagram (\ref{eq:30})) and $pr_{k}\circ \partial=\partial\circ pr_{k-1}$ (Proposition \ref{p:05}). Thus we deduce that $\psi=\bF_{k,\partial}\otimes L^*$ and complete the proof.
\end{proof}

\subsection{Projective geometry}$ $
\medskip

In this subsection, we discuss how to use local sections to compute fundamental forms. This is essentially localizing Gauss maps in the previous subsection at a point. However, we give an alternative elementary and independent approach without quoting Gauss maps. Let us focus on a nonsingular point $x$ and replace $X$ by an affine neighborhood. Recall for a vector field $\partial \in T_X$, it gives rise to a derivation $\partial: V^*\otimes \sO_X\rightarrow V^*\otimes \sO_X$ (see Section 2.1). For a section $B\in V^*\otimes \sO_X$, evaluating $\partial (B)$ at the point $x$ gives a vector $\partial(B)(x)\in V^*$.
\begin{definition} For $u\in T_x(X)$ and $B\in V^*\otimes \sO_X$, we 
	define a vector $\frac{dB}{du}=\partial (B)(x)\in V^*$,
	where $\partial\in T_X$ is a vector field whose evaluation at $x$ is $u$, i.e. $\partial (x)=u$. 
\end{definition}
\begin{lemma} \label{lm:31} The vector $\frac{dB}{du}$ is independent on the choice of $\partial$.
\end{lemma}
\begin{proof} Let $\mf{m}_x$ be the maximal ideal of $\sO_X$ defining the point $x$. Let $\partial'\in T_X$ with $\partial'(x)=u$. Then $\partial-\partial'\in \mf{m}_xT_X$. As $T_X$ has a basis $\{D_1,\ldots,D_n\}$. So we have $\partial-\partial'=\sum a_iD_i$ with $a_i\in \mf{m}_x$. Thus 
	$\partial(B)-\partial'(B)=\sum a_iD_i(B)$. As a consequence, $\partial(B)(x)-\partial'(B)(x)=\sum a_i(x)D_i(B)(x)=0$, i.e., $\partial(B)(x)=\partial'(B)(x)$.
	
\end{proof}

Consider the subsheaf $P^*_{k-1}\otimes L^*$ of $V^*\otimes \sO_X$. By definition, $\widehat{\bT}_x^{k-1}(X)=P^*_{k-1}\otimes L^*\otimes k(x)$ for $x\in \bU_{k-1}$. So there is an evaluation map $P^*_{k-1}\otimes L^*\rightarrow \widehat{\bT}_x^{k-1}(X)$.
\begin{definition} For $x\in \bU_{k-1}$, $v\in \widehat{\bT}^{k-1}_x(X)$ and  $u\in T_x(X)$, let $B_v\in P^*_{k-1}\otimes L^*$ such that  $B_v(x)=v$ and let $\partial_u\in T_X$ such that $\partial_u(x)=u$. Define $\frac{dv}{du}$ to be the image of the vector $\partial_u(B_v)(x)$ in the quotient space $V^*/\widehat{\bT}^{k-1}_x(X)$, i.e.,
	$\frac{dv}{du}:= \Big(\partial_u(B_v)(x) \mod \widehat{\bT}^{k-1}_x(X)\Big)\in V^*/\widehat{\bT}^{k-1}_x(X)$.
\end{definition}
\begin{lemma} $\frac{dv}{du}$ is independent on the choice of $B_v$ and $\partial_u$.
\end{lemma}
\begin{proof} Let $B'_v\in P^*_{k-1}\otimes L^*$ such that $B'_v(x)=v$ and let $\partial'_u\in T_X$ such that $\partial'_u(x)=u$. It suffices to show $\partial_u(B_v)(x)=\partial'_u(B'_v)(x)$ modulo $\widehat{\bT}^{k-1}_x(X)$. But by Lemma \ref{lm:31}, we know that $\partial_u(B'_v)(x)=\partial'_u(B'_v)(x)$.	 Thus it suffices to show 
	$\displaystyle \Big(\partial_u(B_v)(x)\mod \widehat{\bT}^{k-1}_x(X)\Big)=\Big(\partial_u(B'_v)(x) \mod \widehat{\bT}^{k-1}_x(X)\Big)$.
	To this end, note that $B_v-B'_v\in \mf{m}_x(P^*_{k-1}\otimes L^*)$ where $\mf{m}_x$ is the maximal ideal of $\sO_X$ defining $x$. So we can write $B_v-B'_v=\sum a_iB_i$ for some $a_i\in \mf{m}_x$ and $B_i\in P^*_{k-1}\otimes L^*$. Thus $\partial_u(B_v)-\partial_u(B'_v)=\sum \partial_u(a_i)B_i+\sum a_i\partial_u(B_i)$ and thus $\partial_u(B_v)(x)-\partial_u(B'_v)(x)=\sum \partial_u(a_i)(x)B_i(x)\in \widehat{\bT}^{k-1}_x(X)$ which proves the lemma.
\end{proof}

\begin{proposition}\label{p:31} For $x\in \bU_{k-1}$, define the morphism 
$$
		f_{k,x}:T_x(X)\times \widehat{\bT}^{k-1}_x(X)\longrightarrow  V^*/\widehat{\bT}^{k-1}_x(X)
$$
by $f_{k,x}(u,v)=\frac{dv}{du}$ for $u\in T_x(X)$ and $v\in \widehat{\bT}^{k-1}_x(X)$.
Then one has
	\begin{enumerate}
		\item $f_{k,x}$ is $\nC$-bilinear.
		\item If $x\in \bU_k$, then the image of $f_{k,x}$ is $\widehat{\bT}^k_x(X)/\widehat{\bT}^{k-1}_x(X)$, which is also the image of $\bF_{k}\otimes L^*\otimes k(x)$.
		\item For $u\in T_x(X)$, the induced map $f_{k,x}(u): \widehat{\bT}^{k-1}_x(X)\rightarrow  V^*/\widehat{\bT}^{k-1}_x(X)$
		equals $dg_{k-1,x}(u)$ where $g_{k-1,x}$ is the Gauss map $g_{k-1}$ at the point $x$.
		\item For $u\in T_x(X)$ and $v\in \widehat{\bT}^{k-1}_x(X)$, if $v\in  \widehat{\bT}^{k-2}_x(X)$ then $f_{k,x}(u,v)=0$.
	\end{enumerate}
\end{proposition}
\begin{proof} The results follow directly from definition and Proposition \ref{p:05}. So we leave the details for the reader to check. 
	
\end{proof}

\begin{remark}\label{rmk:32} Let $x\in \bU_{k}$. Proposition \ref{p:31} (4) implies that $f_{k,x}$ induces a well-defined $\nC$-bilinear surjective morphism 
	$$\bar{f}_{k,x}:T_x(X)\times \frac{\widehat{\bT}^{k-1}_x(X)}{\widehat{\bT}_x^{k-2}(X)}\longrightarrow \frac{\widehat{\bT}^k_x(X) }{\widehat{\bT}^{k-1}_x(X)}.$$
	Iterating this morphism with $\bar{f}_{k-1,x}$, $\ldots$, $\bar{f}_{1,x}$, one obtains a surjective multi-linear morphism  
	$\gamma_{k,x}:T_x(X)\times \cdots \times T_x(X)\times \hat{x}\rightarrow \frac{\widehat{\bT}^k_x(X)}{\widehat{\bT}^{k-1}_x(X)}$
	where note that $\hat{x}=\widehat{\bT}^0_x(X)$. Equivalently, it determines a map $\gamma'_{k,x}:T_x(X)\times \cdots \times T_x(X)\rightarrow \Hom(\hat{x},\frac{\widehat{\bT}^k_x(X) }{\widehat{\bT}^{k-1}_x(X)})$.
	Recall the twisted fundamental form (Definition \ref{def:02}) $\bF_{k,x}\otimes L^*:S^kT_x(X)\otimes \hat{x}\rightarrow R^*_{k-1,x}(X)\otimes \hat{x}=\frac{V^*}{\widehat{\bT}^{k-1}_x(X)}$ has the image $\widehat{\bT}^k_x(X) /\widehat{\bT}^{k-1}_x(X)$ and induces the fundamental form  
	$\bF_{k,x}: S^kT_x(X)\rightarrow R^*_{k-1,x}(X)=\Hom(\hat{x},\frac{V^*}{\widehat{\bT}^{k-1}_x(X)})$.
	One checks iteratively that $\gamma_{k,x}$ factors through $\bF_{k,x}\otimes L^*$ in a natural way. Thus we obtain the following commutative diagram 
	$$\xymatrix{
		T_x(X)\times \cdots \times T_x(X) \ar[d]\ar[r]^-{\gamma'_{k,x}}&\Hom(\hat{x}, \frac{\widehat{\bT}^k_x(X) }{\widehat{\bT}^{k-1}_x(X)}) \ar@{^{(}_->}[d] \\
		S^{k}T_X   \ar[r]^-{\bF_{k,x}}&	 \Hom(\hat{x},\frac{V^*}{\widehat{\bT}^{k-1}_x(X)}) =R^*_{k-1,x}(X)
	}$$
	For vectors $u_1,\ldots,u_k\in T_x(X)$, we may also write $\bF_{k,x}(u_1,u_2,\ldots, u_k)$ for $\bF_{k,x}(u_1u_2\ldots u_k)$. In addition, $\bF_{k,x}(u_1,u_2,\ldots, u_k)$ is a linear map from $\hat{x}$ to $\frac{V^*}{\widehat{\bT}^{k-1}_x(X)}$ and its image $\bF_{k,x}(v_1,v_2,\ldots, v_k)(\hat{x})$ is either zero or a one-dimensional space in $\frac{V^*}{\widehat{\bT}^{k-1}_x(X)}$.

\end{remark}

\begin{remark}\label{rmk:05} We can  use local sections to compute  $\bF_{k,x}:S^kT_x(X)\rightarrow R_{k-1,x}^*(X)$ at a point $x\in \bU_{k}$. Under the identification $R_{k-1,x}^*(X)=\Hom (\hat{x}, V^*/\widehat{\bT}^{k-1}_x(X))$, for any form $w\in S^kT_x(X)$, $\bF_{k,x}(w)$ is a linear map from $\hat{x}$ to $V^*/\widehat{\bT}^{k-1}_x(X)$. Sub-spaces of the form $\bF_{k,x}(w)(\hat{x})$ span the space $\widehat{\bT}^k_x(X)/\widehat{\bT}^{k-1}_x(X)$. As $\sO_X(-1)$ is a subsheaf of $V^*\otimes \sO_X$, so locally around $x$, we can trivialize $\sO_X(-1)=\sO_XB_0$ by a section $B_0\in V^*\otimes\sO_X$. Note that $\hat{x}=\nC\cdot B_0(x)$. We shall utilize $B_0$ to compute fundamental forms. To be precise,  for $u_1,\ldots, u_k\in T_x(X)$, extending $u_i$ to a vector field $\partial_{u_i}\in T_X$ such that $\partial_{u_i}(x)=u_i$, we define
	$$\frac{d^kB_0}{du_1du_2\ldots du_k}:=(\partial_{u_k}\circ \partial_{u_{k-1}}\circ\ldots \circ \partial_{u_1})(B_0)(x),$$
i.e., apply $\partial_{u_i}$ consecutively to $B_0$ and then evaluate the resulting section at $x$. It is a vector in $V^*$ and it does depends on the choice and the ordering of $\partial_{u_i}$. But inductively by the remark above it is easy but tedious to check that 
$$\bF_{k,x}(u_1\ldots u_k)(B_0(x))=\Big(\frac{d^kB_0}{u_1u_2\ldots u_k}\mod \widehat{\bT}^{k-1}_x(X)\Big)\in V^*/\widehat{\bT}^{k-1}_x(X).$$
	
	For instance, if $x\in \bU_1$, the first fundamental form $\bF_{1,x}:T_x(X)\rightarrow \Hom(\hat{x},V^*/\hat{x})$ is an injective map whose image is $\Hom(\hat{x},\widehat{\bT}_x(X)/\hat{x})$. So we obtain an identification $T_x(X)=\Hom(\hat{x},\widehat{\bT}_x(X)/\hat{x})$ and a vector $u\in T_x(X)$ can be considered as a linear map $u:\hat{x} \rightarrow  V^*/\hat{x}$ defined by $u(B_0(x))=\frac{dB_0}{du} \mod \hat{x}$.  Note that $u(\hat{x})$ is a one-dimensional subspace in $\widehat{\bT}_x(X)/\hat{x}$ and the one-dimensional subspace $\nC\cdot  \frac{dB_0}{du} $ of $\widehat{\bT}_x(X)$ can be thought of as a geometric realization of the vector $u$. In this way, if $\{u_1,\ldots, u_n\}$ is a basis of $T_x(X)$, then we see that $\widehat{\bT}_x(X)=\langle B_0(x), \frac{dB_0}{du_1},\ldots, \frac{dB_0}{du_n}\rangle$.
	Sometimes one can  choose an isomorphism $\hat{x}\cong \nC$ to have an isomorphism $T_x(X)\cong\widehat{\bT}_x(X)/\hat{x}$, which is not canonical but  unique up to scalars.
\end{remark}
\begin{remark} The commutative diagram (\ref{eq:81}) induces the following diagram, which is useful to compute fundamental forms. Let $x\in \bU_{k}$ and let $u\in T_x(X)$ be a vector, then one has 
\begin{equation}\label{eq:91}
	\begin{CD} \displaystyle
		S^{k-1}T_x(X)@>\bF_{k-1,x}>> \displaystyle \Hom(\hat{x},\frac{\widehat{\bT}^{k-1}_x(X)}{\widehat{\bT}^{k-2}_x(X)}) @.\subseteq R^*_{k-2,x}(X)\\
		@V u VV @V u VV\\
		S^kT_x(X)@>\bF_{k,x}>> \displaystyle \Hom(\hat{x},\frac{\widehat{\bT}^{k}_x(X)}{\widehat{\bT}^{k-1}_x(X)})@.\subseteq R^*_{k-1,x}(X),
	\end{CD}
\end{equation}
where the right-hand side $u$ is a naturally induced map. So for vectors $v_1,\ldots, v_{k-1}\in T_x(X)$, we obtain $\bF_{k,x}(u,v_1,\ldots,v_{k-1})=u\circ \bF_{k-1,x}(v_1,\ldots, v_{k-1})$.
	
\end{remark}

\section{Vanishing of fundamental forms}

\noindent In this section, after studying the higher order tangent varieties we give the proof for the vanishing of fundamental forms.  Throughout this section,  $X\subseteq\nP^r$ is a quasi-projective variety of dimension $n\geq 1$ and $V=H^0(\nP^r,\sO_{ \nP^{r}}(1))$.

\subsection{Higher order tangent varieties}

\begin{definition} Let $X\subseteq \nP^r$ be a quasi-projective variety. The {\em $k$-th tangent variety}  $\tau^k(X)$ is defined to be the Zariski closure of the union of projective $k$-th tangent spaces of $X$ at points of $\bU_k$. In particular, the first tangent variety is the tangent variety $\tau(X)$.
\end{definition}

Recall that at a point $x\in \bU_k$ we have a bilinear map $f_{k+1,x}$ defined in Proposition \ref{p:31}. Take a vector $u\in T_x(X)$, it induces a linear map $f_{k+1,x}(u):\widehat{\bT}^k_x(X)\rightarrow V^*/\widehat{\bT}^k_x(X)$. On the other hand, for a point $y\in \bT^k_x(X)$, $\hat{y}$ is a one dimensional subspace of $\widehat{\bT}^k_x(X)$ and we obtain an induced linear map 
\begin{eqnarray*}
	f_{k+1,x,y}:T_x(X) &\longrightarrow &\Hom(\hat{y},V^*/\widehat{\bT}^k_x(X)) \\
	u&\mapsto & f_{k+1,x}(u)|_{\hat{y}}.
\end{eqnarray*}
Observe that geometrically the rank of the map $f_{k+1,x,y}$ is the same as the dimension of the linear space in $V^*/\widehat{\bT}^k_x(X)$ spanned by the image $f_{k+1,x}(u)(\hat{y})$ for all $u\in T_x(X)$. In addition, the rank of $f_{k+1,x,y}$ has an obvious upper bound $\dim X$.

\begin{theorem}\label{thm:52} Let $X\subseteq \nP^r$ be a quasi-projective variety and let $k\geq 1$. There exists an open subset $\bU_{\tau^k}$ contained in $\bU_k$ satisfying the following property. For each point $x\in \bU_{\tau^k}$, there exists an open set $U_x$ of $\bT^k_x(X)$ such that for a point $y\in U_x$, $\tau^k(X)$ is nonsingular at $y$ and 
	$$\uIm f_{k+1,x,y}= \frac{T_y(\tau^k(X)) }{T_y(\bT^k_x(X))}.$$
 In particular, 
	$$\dim \tau^k(X)=\dim \bT^k_x(X)+\rank f_{k+1,x,y}.$$

\end{theorem}
\begin{proof} We replace $X$ by the open set $\bU_k$. Then $P_k(1)$ is a locally free sheaf on $X$ of rank $\mf{t}_k+1$. Let $\nG=\nG(\nP^r, \mf{t}_k)$ be the Grassmannian parameterizing  $\mf{t}_k$-dimensional linear spaces in $\nP^r$. Let $I$ be the universal family over $\nG$. Write $M=\nP(P_k(1))$. One has a commutative diagram 
\begin{equation}
	\begin{CD}
		\xymatrix{
			M\ar[d]_{\pi}\ar[r]^{g'}  \ar@/^1.2pc/[rr]|{\tau_k}&	I \ar[d] \ar[r]^{p} &\nP^r\\
			X\ar[r]^{g_k}&	\nG & 
		}
	\end{CD}
\end{equation}
	where $g_k$ is the $k$-th Gauss map and $M=X\times_{\nG}I$. The composition map $\tau_k:=p\circ g'$ is given by the tautological bundle on $M$ and the closure of its image is the variety $\tau^k(X)$. 
	
	By the generic smoothness, there exists an open subset $U_M$ of $M$ on which the map $\tau_k$ is smooth and $\tau_k(U_M)$ is in the nonsingular locus of $\tau^k(X)$. We take $U_{\tau^k}=\pi(U_M)$ and for each $x\in U_{\tau^k}$ we take $U_x=M_x\cap U_M$ where $M_x$ is the fiber over $x$ and is the projective tangent space $\bT^k_x(X)$.	Let $z=(x,y)\in M$ be a point such that $x\in U_{\tau^k}$ and $y\in U_x\subset M_x=\bT^k_x(X)$. To save notations, write $L=\bT^k_x(X)$. Let $z'=g'(z)$ and consider $y=p(z')\in \nP^r$. By construction $y$ is a nonsingular point of $\tau^k(X)$. The Gauss map $g_k$ sends $x$ to the point $[L]\in \nG$. One has the following diagram of the Zariski tangent spaces 
	\begin{equation}\label{eq:51}
		\begin{CD}
		\xymatrix{
				T_y(L)\ar@{=}[r] \ar@{^{(}->}[d] & T_y(L)\ar@{=}[r]\ar@{^{(}->}[d] &T_y(L)\ar@{^{(}->}[d] \\
				T_z(M)\ar@{->>}[d]\ar[r]^{dg'} &	T_{z'}(I) \ar@{->>}[d] \ar[r]^{dp} &T_y(\nP^r)\ar@{->>}[d]\\
				T_x(X)\ar[r]^{dg_k}&	T_{[L]}(\nG)\ar[r]^{e} & N_y(L) 
			}
		\end{CD}
	\end{equation}
	where $e$ is an induced map. Note that $d\tau_k=dp\circ dg'$ and the Snake Lemma shows that 
	$\ker d\tau_k=\ker (e\circ dg_k)$. Recall that 
	$T_{[L]}(\nG)=H^0(N_L)$ where $N_L$ is the normal sheaf of $L$ in $\nP^r$. The map $e$
	turns out to the evaluation map of global sections of $N_L$ at the point $y$. Since  $N_L=(V^*/\hat{L})\otimes \sO_L(1)$,  $H^0(N_L)=\Hom(\hat{L}, V^*/\hat{L})$, and $N_y(L)=\Hom(\hat{y}, V^*/\hat{L})$, the map $e$ is given by 
	\begin{eqnarray*}
		e: \Hom(\hat{L},V^*/\hat{L}) &\longrightarrow & \Hom(\hat{y},V^*/\hat{L})\\
		\varphi&\mapsto & \varphi|_{\hat{y}}.
	\end{eqnarray*}
	So for $u\in T_x(X)$, $e\circ dg_k(u)=f_{k+1,x}(u)|_{\hat{y}}$ and therefore
	$e\circ dg_k=f_{k+1,x,y}$. As a consequence 
	$\ker\tau_k=\ker f_{k+1,x,y}$.
	Since the image of $d\tau_k$ is $T_y(\tau^k(X))$, the right-hand-side vertical sequence in (\ref{eq:51}) gives that $T_y(\tau^k(X))/T_y(L)=\im (e\circ dg_k)=\im f_{k+1,x,y}$.
	This finishes the proof.
\end{proof}

\begin{remark}\label{rmk:41}
	The expected dimension of $k$-th tangent variety $\tau^k(X)$ is $\dim \bT^{k}_x(X)+n$ where $x$ is a generic point of $X$. We say $\tau^k(X)$ is {\em degenerate} if its dimension is smaller than the expected dimension (in this definition, it is possible that $\tau^k(X)$ fills up the whole space $\nP^r$ if $r$ is smaller than the expected dimension of $\tau^k(X)$). We can define the {\em $k$-th tangent defect} to be the number $\rank P_k+\dim X-\dim \tau^k(X)-1$.
	In particular, the tangent defect is the number $2\dim X-\dim \tau(X)$. Except of the projective tangent space $\bT_x(X)$, it is hard in general to calculate the dimension of $\bT^k_x(X)$. However, by Theorem \ref{thm:52}, we see that the $k$-th tangent defect equals $\dim X-\rank f_{k+1,x,y}$. Thus the fundamental form can be used to detect the degeneracy of $\tau^k(X)$, as described in the following corollary.
\end{remark}

\begin{corollary} Let $X\subseteq\nP^r$ be a quasi-projective variety and $x\in X$ be a generic point. If
	$\rank \bF_{k+1,x}<\dim  X$
	then the $k$-th tangent variety $\tau^k(X)$ is degenerate. 	
\end{corollary}
\begin{proof} We use  notations in Theorem \ref{thm:52}. Recall that $\bF_{k+1,x}$ is induced by the twisted fundamental form $\bF_{k+1,x}\otimes L^*:S^{k+1}T_x(X)\otimes \hat{x}\rightarrow V^*/\widehat{\bT}^k_x(X)$. Thus $\rank \bF_{k+1,x}=\rank \bF_{k+1,x}\otimes L^*=\dim \bF_{k+1,x}\otimes L^*(S^{k+1}T_x(X)\otimes \hat{x})$.
	By Remark \ref{rmk:32}, there is a commutative diagram 
	$$\xymatrix{
	T_x(X)\times \cdots \times T_x(X)\times \hat{x}\ar[r]\ar[d]	&\displaystyle T_x(X)\times \frac{\widehat{\bT}^k_x(X)}{\widehat{\bT}^{k-1}_x(X)}\ar[d]^-{\bar{f}_{k+1,x}}\\
	S^{k+1}T_x(X)\otimes \hat{x}	 \ar[r]_-{\bF_{k+1,x}\otimes L^*}	&	\displaystyle \frac{V^*}{\widehat{\bT}^k_x(X)}\\
	}$$
Observe that the rank of $f_{k+1,x,y}$ is the dimension of the space $f_{k+1,x}(T_x(X)\times \hat{y})$ which is the same as the space $\bar{f}_{k+1,x}(T_x(X)\times ( \hat{y} \mod \widehat{\bT}^{k-1}_x(X)))$. But the latter one is contained in the space $\bF_{k+1,x}\otimes L^*(S^{k+1}T_x(X)\otimes \hat{x})$.	Hence we conclude that $\rank f_{k+1,x,y}\leq \rank \bF_{k+1,x}<\dim X$.
	Then by Theorem \ref{thm:52}, this means $\dim \tau^k(X)< \dim \widehat{\bT}^k_x(X)+n$ and thus $\tau^k(X)$ is degenerate. 
\end{proof}

In the special case of the theorem when $k=1$, the map $f_{2,x}$ factors through $\bar{f}_{k,2}:T_x(X)\times \widehat{\bT}_x(X)/\hat{x}$. Composed with $f_{1,x}$,  we see that there is a vector $w\in T_x(X)$ such that $f_{1,x}(w)(\hat{x})=(\hat{y}\mod \hat{x})$. Equivalently, if we choose an isomorphism $\hat{x}\cong \nC$ so that $T_x(X)\cong \widehat{\bT}_x(X)/\hat{x}$, then $w$ is in the class of $\hat{y}$. Geometrically, one can view $w$ as a vector in the direction determined by the line $\overline{xy}$ connecting the points $x$ and $y$. In this way, we see 
$\rank f_{2,x,y}=\rank \bII_{x,w}$.
Then we obtain the following corollary about the tangent variety in terms of the second fundamental form.
\begin{corollary}\label{p:09} Let $X\subseteq \nP^r$ be a quasi-projective variety. There exists an open subset $\bU_{\tau}$ contained in $\bU_1$ of $X$ satisfying the following property. For each point $x\in \bU_{\tau}$, there exists an open set $U_x$ of the projective tangent space $\bT_x(X)$ such that for a point $y\in U_x$ let $w\in T_x(X)$ be a nonzero vector in the direction of the line $\overline{xy}$, then 
	$$\uIm \bII_{x,w}=\frac{T_y(\tau(X))}{T_y(\bT_x(X))},$$
	where $\bII_{x,w}=\bII_x(\_,w):T_x(X)\rightarrow N_x(X)$ induced by the second fundamental form $\bII_x$ at $x$. In particular, 
	$$\dim \tau(X)=\dim X+\rank \bII_{x,w}.$$
\end{corollary}
\begin{remark} The result in corollary was proved by Griffiths-Harris \cite[5.5]{Griffiths:AGandLDG} under the condition $r\geq 2n$. It was also proved by Landsberg  \cite[4.10]{Landsberg:DegSecTan} using his notion of $II$-generic vector.
\end{remark}

\begin{corollary} Let $X\subseteq\nP^r$ be a quasi-projective variety.  If
	$\rank \bII_x<  \dim X$	at a generic point $x$ then the tangent variety $\tau(X)$ is degenerate. 	
\end{corollary}
\begin{proof} As $\rank \bII_{x,w}\leq \rank \bII_x$ for any vector $w\in T_x(X)$, the result follows from the theorem.	
\end{proof}

\bigskip
\subsection{Vanishing of fundamental forms}$ $
\medskip

\noindent In this subsection, we prove a vanishing theorem for fundamental forms. 

\begin{lemma}[{\bf Curve Selection}]\label{p:10} Let $X\subseteq\nP^r$ be a quasi-projective nonsingular variety. Suppose that $Z$ is a curvilinear  subscheme of finite length supported at a point $p\in X$. For any point $q\in X$, there exists a nonsingular irreducible curve $C$ in $X$ such that $C$ contains both $Z$ and $q$. In particular, if $U\subseteq X$ is an open subset and $x\in X$ is a closed point, then for any nonzero vector $u\in T_x(X)$, there exists a nonsingular curve $C$ in $X$ passing through $x$ such that $u\in T_x(C)\subseteq T_x(X)$ and $C\cap U\neq \emptyset$.
\end{lemma}
\begin{proof} By embedding $X$ in a projective variety and then resolving singularities, we may assume $X$ is projective. By induction on the dimension of $X$, it suffices to show that there exists an nonsingular irreducible hypersurface of $X$ contains both $Z$ and $q$. To construct such hypersurface, note first that if the length of $Z$ is one, then it follows from the Bertini theorem by blowing up both $p$ and $q$ and using a sufficient positive very ample line bundle. So in the sequel, we assume $\length Z\geq 2$.  Let $\pi:\widetilde{X}\rightarrow X$ be the blowup of $X$ along the point $p$ with the exceptional divisor $E$. Since $Z$ is curvilinear, a local calculation shows that $I_Z\cdot\sO_{\widetilde{X}}=I_{\widetilde{Z}}(-E)$, where $I_Z$ is the defining ideal of $Z$ and $I_{\widetilde{Z}}$ is an ideal sheaf of $\sO_{\widetilde{X}}$ defining a curvilinear subscheme of length $\length Z-1$ supported at a point on $E$. Let $L$ be a very ample line bundle on $X$ such that $\pi^*L(-E)$ is also very ample. Then iteratively, we can have a general nonsingular irreducible hyperplane $\widetilde{H}\in |\pi^*(-E)|$ containing both $\widetilde{Z}$ and $q$. Then let $H=\pi(\widetilde{H})$ which is a general nonsingular hypersurface passing through the point $p$ and $q$. It is clear that $I_H\cdot \sO_{\widetilde{X}}=I_{\widetilde{H}}(-E)\subseteq I_{\widetilde{Z}}(-E)=I_Z\cdot \sO_{\widetilde{X}}$. Since $Z$ is curvilinear one checks that   $\pi_*(I_Z\cdot \sO_{\widetilde{X}})=I_Z$ and thus $I_H\subseteq I_Z$ as desired.
	
\end{proof}

\begin{lemma}\label{p:07} Consider $x\in C\subseteq X$, where $X$ is a nonsingular affine variety, $x$ is a closed point and $C$ is a nonsingular curve passing through $x$. Let $u\in T_x(C)\subseteq T_x(X)$ be a nonzero tangent vector and let $\partial_u\in T_X$ be a vector field such that $\partial_{u}(x)=u$. Shrinking $X$ if necessary, there exists an element $t\in \sO_C$ such that $t$ generates the maximal ideal of $x$ in $\sO_C$, $dt$ generates $\Omega^1_C$, and the dual $\partial_t$ of $dt$ is the restriction of $\partial_u$, i.e. $\partial_t=\partial_u\otimes 1$ under the inclusion $T_C\hookrightarrow T_X|_C$.
\end{lemma}
\begin{proof} In the short exact sequence $0\rightarrow T_x(C)\rightarrow T_x(X)\rightarrow N_x(C)\rightarrow 0$, the image of $u$ in $N_x(C)$ is zero. Thus we can shrink $X$ if necessary such that the restriction of $\partial_u$ on to $C$ is in $T_C$ under the inclusion $T_C\hookrightarrow T_X|_C$. The vector $u$ is a base for $T_x(C)$. Write $u^*$ the dual base for $\Omega^1_x(C)=\mf{m}/\mf{m}^2$ where $\mf{m}$ is the maximal ideal of $\sO_C$ defining $x$. We can take an element $t\in \mf{m}$ such that its image in $\Omega^1_x(C)$ is $u^*$. By shrinking $X$ again, we assume $t$ generates the maximal ideal $\mf{m}$ in $\sO_C$ and $dt$ generate $\Omega^1_C$, and the dual $\partial_t$ of $dt$ is the restriction of $\partial_u$, i.e. $\partial_t=\partial_u\otimes 1$.
	
\end{proof}
\begin{remark}\label{rmk:03} In the above proposition,  $\partial_t$ and $\partial_u$ give the same tangent vector $u\in T_x(C)\subseteq T_x(X)$ at the point $x$. 	Let $a\in \sO_X$ be a section and write $\bar{a}=a\otimes 1\in \sO_C$ the restriction of $a$ onto $C$. We have 
	$\overline{\partial_u(a)}=\partial_t(\bar{a})\in \sO_C$. In iterative way we obtain 
	$\overline{\partial^m_u(a)}=\partial^m_t(\bar{a})$, for $m\geq 0$.
	The  evaluations of above sections at $x$ is denoted by
	$$\frac{d^ma}{du^m}:=\partial^m_u(a)(x), \text{ and }\frac{d^m\bar{a}}{dt^m}:=\partial^m_t(\bar{a})(x).$$
	So we conclude that $\frac{d^ma}{du^m}=\frac{d^m\bar{a}}{dt^m}$ for all $m\geq0$.
\end{remark}

\begin{definition}Let $X\subseteq \nP^r$ be a quasi-projective variety. For $k\geq 1$, define the number 
	$$\delta_k=\dim \bT^k_x(X)\cap \bT^k_y(X), \text{ for generic points }x, y\in X,$$
	where we use convention that $\delta_k=-1$ if $\bT^k_x(X)\cap \bT^k_y(X)=\emptyset$.
	
\end{definition}

\begin{proposition} Let $X\subseteq\nP^r$ be a quasi-projective variety and let $k\geq 1$. There exists an open subset $U$ contained in $\bU_k$ such that for $x\in U$, there exists an open subset $U_x\subseteq U$ with the property that 
	$$\delta_k=\dim \bT^k_x(X)\cap \bT^k_y(X), \text{ for any }y\in U_x.$$
\end{proposition}
\begin{proof} Without loss of generality, we replace $X$ by the open subset $\bU_k$. Consider $Y=X\times X$ with the natural projections $p_1$ and $p_2$ to $X$. Then one has the surjective morphism 
	$$V\otimes \sO_Y\longrightarrow p_i^*P_k(1)$$
	induced by pulling back the Taylor series map $V\otimes \sO_X\rightarrow P_k(1)$. So the variety $\nP(p_i^*P_k(1))$, $i=1$ and $2$, is a subvariety of $\nP(V\otimes \sO_Y)$. Let $Z=\nP(p_1^*P_k(1))\cap \nP(p_2^*P_k(1))$ with the projection $\pi:Z\rightarrow Y$. Then there exists an open subset $U_Y$ of $Y$ such that the fiber of the induced map $Z_U\rightarrow U_Y$ is a linear space of dimension $\delta_k$, where $Z_U=\pi^{-1}(U_Y)$. We take $U=p_1(U_Y)$ and for each $x\in U$, we take $U_x=U_Y\cap \{x\}\times X$. This proves the proposition.  
\end{proof}

\begin{remark} The case $k=1$ of the proposition is essentially the Terracini lemma on secant varieties. Recall that the secant defect $\delta_X:=2\dim X+1-\dim \Sec(X)$.
The secant variety is degenerate if $\delta_X>0$. Directly by definition, we have $\delta_1=\delta_X-1$. 
\end{remark}

\begin{definition}Let $X\subseteq \nP^r$ be a quasi-projective variety. For $k\geq 1$, we define a number 
	$$\theta_k=\rank f_{k+1,x}(u)$$
where $x\in X$ is a generic point, $u\in T_x(X)$ is a generic vector and $f_{k+1,x}(u): \widehat{\bT}^{k}_x(X)\rightarrow V^*/\widehat{\bT}^k_x(X)$ is the map induced by $f_{k+1,x}$. Observe that $\theta_k\leq \mf{t}_k-\mf{t}_{k-1}$.
\end{definition}

\begin{remark}\label{rmk:53} The number $\theta_k$ is a well-defined number at generic points. Indeed, replace $X$ by an affine open subset in $\bU_k$. Consider the space $T=\Spec \Sym(T^*_X)$ and $H=\Spec(\Sym \sHom(P^*_k(-1),R^*_k(-1))^*)$ over $X$. The morphism $dg_k$ induces a morphism $g:T\rightarrow H$  over $X$ such that
for a point $x\in X$, the fiber $T_x=T_x(X)$, $H_x=\Hom(P^*_x(X),R^*_x(X))$ and the morphism $g_x=dg_{k,x}$. A closed point $h\in H$ is a morphism in $H_x$ where $x=q(h)$ and $\rank h$ is the rank of $h$ as  a linear map. Thus for each $i$ we can define a closed subset $H_i=\{h\in H\mid \rank h\leq i\}$ of $H$ such that $H_0\subseteq H_1\subseteq \ldots H$. So there exists $t$ such that $g(T)\subseteq H_t$ but $g(T)\nsubseteq H_{t-1}$. It is clear that the number $\theta_k=t$. Let  $U_T=g^{-1}(H_t-H_{t-1})$ which is  an open subset of $T$ and let $U=p(U_T)$ be the open subset in $X$. Then we see that for any $x\in U$, we can obtain an open subset $U_x=U_T\cap T_x$ in the Zariski tangent space $T_x(X)$ such that for any  $u\in U_x$, $f_{k+1,x}(u)=g(u)$ has rank $\theta_k$. Furthermore, if we pick a vector field $\partial\in T_X$ such that $\partial(x)\in U_x$, then for any point $x'$ in a neighborhood of $x$, $\partial(x')\in U_{x'}$.
	
\end{remark}

%
\begin{proof}[{\bf Proof of Theorem \ref{thm:vanishing}}] Working on an affine open subset in $\bU_k$ we may assume $X=\Spec \sO_X$ is affine.  Let $x\in X$ be a point.  Take a generic vector $u\in T_x(X)$ such that $\theta_k=\rank f_{k+1,x}(u)$.  By Curve Selection Lemma \ref{p:10} there exits a nonsingular curve $C\subseteq X$ passing through $x$ such that $u\in T_x(C)$ and 
	\begin{equation}\label{p:13}
		\delta_k=\dim \bT^k_x(X)\cap\bT^k_p(X), \text{ for }p\in C-\{x\}.
	\end{equation}

	Take a vector field $\partial_u\in T_X$ such that $\partial_{u}(x)=u$. By Proposition \ref{p:07},  there exists an element $t\in \sO_C$ such that $t$ generates the maximal ideal of $x$ in $\sO_C$, $dt$ generates $\Omega^1_C$, and the dual $\partial_t$ of $dt$ is the restriction $\partial_u\otimes 1$.
	
	Recall for a section $B\in V^*\otimes \sO_X$ we use notation $\displaystyle\frac{d^mB}{du^m}:=\partial^m_uB(x)$, for $m\geq 0$. Denote by $\bar{B}=B\otimes 1\in V^*\otimes \sO_C$  the restriction of $B$ onto $C$, and write 
	$\displaystyle \frac{d^m\bar{B}}{dt^m}:=\partial^m_t\bar{B}(x)$ for $m\geq 0$.

	Now we consider free modules $\sO_X(-1)\subseteq P^*_k(-1)\subseteq V^*\otimes \sO_X$. There exists sections $B_i$ of $V^*\otimes\sO_X$, $i=0,\ldots, \mf{t}_k$, such that 
	$$\sO_X(-1)=\sO_XB_0, \text{ and }P^*_k(-1)=\sO_XB_0\oplus \sO_XB_1\oplus\ldots\oplus\sO_XB_{\mf{t}_k}.$$
	Write $A_0=B_0(x),\ldots, A_{\mf{t}_k}=B_{\mf{t}_k}(x)$ considered as vectors in $V^*$ and note that
	$$\widehat{\bT}^k_x(X)=\langle A_0,\ldots, A_{\mf{t}_k}\rangle.$$
	Write $\widehat{\bT}^k_x(X)+\uIm f_{k+1,x}(u)$ as the subspace of $V^*$ over the image $f_{k+1,x}(u)(\widehat{\bT}^{k}_x(X))$. Then 
	$$\widehat{\bT}^k_x(X)+\uIm f_{k+1,x}(u)=\langle A_0,\ldots, A_{\mf{t}_k}, \frac{dB_1}{du},\ldots, \frac{dB_{\mf{t}_k}}{du}\rangle\subseteq \widehat{\bT}^{k+1}_x(X).$$
	Since $\dim (\widehat{\bT}^k_x(X)+\uIm f_{k+1,x}(u))=\mf{t}_k+1+\mf{t}_k-\delta_k$ by assumption, without loss of generality, we may assume the vectors
	\begin{equation}\label{p:15}
		A_0,\ \ldots,\  A_{\mf{t}_k},\ \frac{dB_1}{du},\ \frac{dB_2}{du},\ \ldots,\ \frac{dB_{\mf{t}_k-\delta_k}}{du}, \text{ are linearly independent}.
	\end{equation}
	
	\begin{claim}\label{p:14} The sections 	$\bar{B}_0,\ \bar{B}_1,\ \ldots,\ \bar{B}_{\mf{t}_k},\ \partial_t\bar{B}_1,\ \partial_t\bar{B}_2,\ \ldots,\ \partial_t\bar{B}_{\mf{t}_k-\delta_k} $
		are linearly independent at each point of $C$ and their span contains sections  $\partial_t\bar{B}_{\mf{t}_k-\delta_k+1},\ \ldots,\ \partial_t\bar{B}_{\mf{t}_k}$. 	
	\end{claim}
	
	{\em Proof of Claim}. To prove the claim, we note that being linearly independent is an open condition and thus (by shrinking $X$ and $C$ if necessary) $\bar{B}_0, \bar{B}_1,\ldots, \bar{B}_{\mf{t}_k}, \partial_t\bar{B}_1,\partial_t\bar{B}_2,\ldots, \partial_t\bar{B}_{\mf{t}_k-\delta_k} $ are linear independent at each point of $C$ because of (\ref{p:15}). On the other hand, the vector field $\partial_u$ (and hence $\partial_t$) gives generic vectors at each point of $C$ (see Remark \ref{rmk:53}). Thus applying Theorem \ref{p:09} at each point $p\in C$, we see that $$\partial_t\bar{B}_{\mf{t}_k-\delta_k+1}(p),\ldots,\partial_t\bar{B}_{\mf{t}_k}(p)\in \langle \bar{B}_0(p), \bar{B}_1(p),\ldots, \bar{B}_{\mf{t}_k}(p), \partial_t\bar{B}_1(p),\partial_t\bar{B}_2(p),\ldots, \partial_t\bar{B}_{\mf{t}_k-\delta_k}(p) \rangle.$$ This proves the claim.
	
	\medskip
	
	We throw away $\delta_k$ sections $B_{\mf{t}_k+1-\delta_k}$, $\ldots$, $B_{\mf{t}_k}$ and construct a section
	$$\sigma:=(A_0\wedge \ldots \wedge A_{\mf{t}_k})\otimes 1\wedge B_0\wedge\ldots \wedge B_{\mf{t}_k-\delta_k}\in (\det W\otimes \sO_X)\wedge (\wedge^{\mf{t}_k+1-\delta_k}V^*\otimes \sO_X),$$ 
	 where $W=\widehat{\bT}^k_x(X)\subseteq V^*$,
	Restricting $\sigma$ onto the curve $C$  yields a section 
	$$\bar{\sigma}=(A_0\wedge \ldots \wedge A_{\mf{t}_k})\otimes 1\wedge \bar{B}_0\wedge\ldots \wedge \bar{B}_{\mf{t}_k-\delta_k}\in (\det W\otimes \sO_C)\wedge (\wedge^{\mf{t}_k+1-\delta_k}V^*\otimes \sO_C).$$
	Since for $p\in C$, $\widehat{\bT}^k_p(X)=\langle \bar{B}_0(p),\ldots, \bar{B}_{\mf{t}_k}(p)\rangle$ intersects $\widehat{\bT}^k_x(X)=\langle A_0,\ldots,A_{\mf{t}_k}\rangle$ in a space of dimension $\delta_k+1$. Thus $\bar{\sigma}(p)=0$ and as a consequence $\bar{\sigma}=0$ on $C$, i.e., 
	\begin{eqnarray}\label{eq:20}
		(A_0\wedge \ldots \wedge A_{\mf{t}_k})\otimes 1\wedge \bar{B}_0\wedge\ldots \wedge \bar{B}_{\mf{t}_k-\delta_k}=0.
	\end{eqnarray}
	Pass to the completion  $\hat{\sO}_{C,x}$ of $\sO_{C,x}$ and take the Taylor expansion for each $\bar{B}_i$,
	\begin{eqnarray*}
		\bar{B}_0&=&A_0+t\frac{d\bar{B}_0}{dt}+\frac{t^2}{2}\frac{d^2\bar{B}_0}{dt^2}+\frac{t^3}{3!}\frac{d^3\bar{B}_0}{dt^3}+\ldots,\\
		\bar{B}_1&=&A_1+t\frac{d\bar{B}_1}{dt}+\frac{t^2}{2}\frac{d^2\bar{B}_1}{dt^2}+\frac{t^3}{3!}\frac{d^3\bar{B}_1}{dt^3}+\ldots, \\
		&\cdots&\\
		\bar{B}_{\mf{t}_k-\delta_k}&=&A_{\mf{t}_k-\delta_k}+t\frac{d\bar{B}_{\mf{t}_k-\delta_k}}{dt}+\frac{t^2}{2}\frac{d^2\bar{B}_{\mf{t}_k-\delta_k}}{dt^2}+\frac{t^3}{3!}\frac{d^3\bar{B}_{\mf{t}_k-\delta_k}}{dt^3}+\ldots.
	\end{eqnarray*}
	Substitute Taylor expansions above into the equation (\ref{eq:20}) and in the result all the coefficients of powers of $t$ are zero. As $\frac{d\bar{B}_0}{dt},\ldots,\frac{d^k\bar{B}_0}{dt^k}\in \langle A_0,\ldots,A_{\mf{t}_k}\rangle$, so in particular, the coefficient of $t^{k+2+\mf{t}_k-\delta_k}$ would be of the form 
	\begin{equation}\label{eq:21}
		\Lambda_0\wedge\Big(\frac{1}{(k+2)!}\frac{d^{k+2}\bar{B}_0}{dt^{k+2}}\wedge \frac{d\bar{B_1}}{dt}\wedge\cdots \wedge \frac{d\bar{B}_{\mf{t}_k-\delta_k}}{dt}+\frac{1}{2(k+1)!}\Delta\Big)=0, \text{ where }
	\end{equation}
    \begin{eqnarray*}
	\Lambda_0&=&A_0\wedge A_1\wedge \cdots\wedge A_{\mf{t}_k}, \text{ and }
\end{eqnarray*}
	$$\Delta=\frac{d^{k+1}\bar{B}_0}{dt^{k+1}}\wedge\frac{d^2\bar{B}_1}{dt^2}\wedge\ldots\wedge\frac{d\bar{B}_{\mf{t}_k-\delta_k}}{dt}+\frac{d^{k+1}\bar{B}_0}{dt^{k+1}}\wedge\frac{d\bar{B}_1}{dt}\wedge \frac{d^2\bar{B}_2}{dt^2}\wedge\ldots \wedge \frac{d\bar{B}_{\mf{t}_k-\delta_k}}{dt}+\ldots+\frac{d^{k+1}\bar{B}_0}{dt^{k+1}}\wedge\frac{\bar{B}_1}{dt}\wedge\ldots\wedge \frac{d^2\bar{B}_{\mf{t}_k-\delta_k}}{dt}.$$ 
	We can calculate $\displaystyle \frac{d^k\bar{B}_0}{dt^k}$, $\displaystyle\frac{d^{k+1}\bar{B}_0}{dt^{k+1}}$, and $\displaystyle\frac{d^{k+2}\bar{B}_0}{dt^{k+2}}$ in terms of $\bar{B}_1,\ldots, \bar{B}_{\mf{t}_k-\delta_k}$. Indeed, notice that  
	$$\partial^k_t\bar{B}_0=\overline{\partial^k_uB_0}\in P^*_k\otimes L^*|_C=\sO_C\bar{B}_0+\ldots+\sO_C\bar{B}_{\mf{t}_k}$$
	and hence we have 
	\begin{eqnarray*}
		\partial^k_t\bar{B}_0&=&a'_0\bar{B}_0+\ldots +a'_{\mf{t}_k}\bar{B}_{\mf{t}_k}, \text{ for some  }a'_i\in \sO_C,\\
		\partial^{k+1}_t\bar{B}_0&=&\sum_{i=0}^{\mf{t}_k} b_i\bar{B}_i+a_1\partial_t\bar{B}_1+\ldots+a_{\mf{t}_k-\delta_k}\partial_t\bar{B}_{\mf{t}_k-\delta_k}, \text{ for some }a_i,b_i\in \sO_C, \text{ (by Claim \ref{p:14})},\\
		\partial^{k+2}_t\bar{B}_0&=&\sum_{i=0}^{\mf{t}_k}c_i\bar{B}_i+\sum_{i=1}^{\mf{t}_k-\delta_k}e_i\partial_t\bar{B}_i+a_1\partial^{2}_t\bar{B}_1+\ldots+a_{\mf{t}_k-\delta_k}\partial^{2}_t\bar{B}_{\mf{t}_k-\delta_k}, \text{ for some }c_i,e_i\in \sO_C.
	\end{eqnarray*}
	Evaluate above sections at $x$, we have
	\begin{eqnarray*}
		\frac{d^k\bar{B}_0}{dt^k}&=& a'_0A_0+\ldots +a'_{\mf{t}_k}A_{\mf{t}_k}, \text{ where }a'_i\in k(x),\\
		\frac{d^{k+1}\bar{B}_0}{dt^{k+1}}&=&\sum_{i=0}^{\mf{t}_k} b_iA_i+a_1\frac{d\bar{B}_1}{dt}+\ldots+a_{\mf{t}_k-\kappa}\frac{d\bar{B}_{\mf{t}_k-\delta_k}}{dt}, \text{ for some }a_i,b_i\in k(x),\\
		\frac{d^{k+2}\bar{B}_0}{dt^{t+2}}&=&\sum_{i=0}^{\mf{t}_k}c_iA_i+\sum_{i=1}^{\mf{t}_k-\delta_k}e_i\frac{d\bar{B}_1}{dt}+a_1\frac{d^{2}\bar{B}_1}{dt^{2}}+\ldots+a_{\mf{t}_k-\delta_k}\frac{d^{2}\bar{B}_{\mf{t}_k-\delta_k}}{dt^{2}}, \text{ for some }c_i,e_i\in k(x).
	\end{eqnarray*}
	Using equality above, we calculate  
	\begin{eqnarray*}
		\Lambda_0\wedge \Delta&=&-\Lambda_0\wedge (a_1\frac{d^{2}\bar{B}_1}{dt^{2}}+\ldots+a_{\mf{t}_k-\delta_k}\frac{d^{2}\bar{B}_{\mf{t}_k-\delta_k}}{dt^{2}})\wedge \frac{d\bar{B}_1}{dt}\wedge\ldots\wedge \frac{d\bar{B}_{\mf{t}_k-\delta_k}}{dt}\\
		&=&-\Lambda_0\wedge \frac{d^{k+2}\bar{B}_0}{dt^{k+2}}\wedge \frac{d\bar{B}_1}{dt}\wedge\ldots\wedge \frac{d\bar{B}_{\mf{t}_k-\delta_k}}{dt}.
	\end{eqnarray*}
	Hence the equation (\ref{eq:21}) becomes 
	$$(\frac{1}{(k+2)!}-\frac{1}{2(k+1)!})\Lambda_0\wedge \frac{d^{k+2}\bar{B}_0}{dt^{k+2}}\wedge \frac{d\bar{B}_1}{dt}\wedge\ldots\wedge \frac{d\bar{B}_{\mf{t}_k-\delta_k}}{dt}=0$$
	Thus 
	$$\frac{d^{k+2}\bar{B}_0}{dt^{k+2}}\in \langle A_0,\ldots, A_{\mf{t}_k}, \frac{d\bar{B}_1}{dt},\ldots, \frac{d\bar{B}_{\mf{t}_k-\delta_k}}{dt}\rangle$$
	This means 
	$$\frac{d^{k+2}\bar{B}_0}{dt^{k+2}}=\frac{d^{k+2}B_0}{du^{k+2}}=\bF_{k+2,x}(u,\ldots,u)=0 \ \mod \widehat{\bT}_x^{k+1}(X).$$
	The result then follows from Lemma \ref{lm:52} and Lemma \ref{lm:53} below.
\end{proof}

\begin{lemma}\label{lm:52} Let $T=\nA^n_{\nC}$ be a vector space considered as an affine space and let $k\geq 1$ be an integer. Let $U\subseteq T$ be an Zariski open subset. Then there exists a basis $\{v_1,\cdots, v_n\}$ of $T$  such that for any integer $0\leq J\leq k$ the sum of $J$ vectors $\sum^J_{j=1}v_{i_j}$ is in $U$.
\end{lemma}
\begin{proof} The open subset $U$ contains a basis for $T$ because, otherwise, its linear span is a proper linear space which is contradict to the density property of Zariski open sets. Let $\{v_1,\ldots,v_n\}$ be a basis contained in $U$. We can re-scale each $\lambda_i v_i$ with a scalar $\lambda_i\in \nC$ to obtain a new basis with the desired property.  This can be done in an iterated way. We start with $v_1$. Since $\nC v_1$ is an affine line so there is only finitely many values of $\lambda_1$ such that $\lambda_1 v_1$ is not in $U$. Hence we can choose a value $\lambda_1$ such that after resetting $v_1=\lambda_1v_1$, $v_1$ has the property that  $v_1, 2v_2,\ldots, kv_1$ are all contained in $U$. Let $B=\{v_1,\ldots, v_t\}$ be the set containing re-scaled vectors. We will re-scale  $v_{t+1}$ and add the resulting vector into $B$. To do this, for any integer $0\leq J\leq k-1$ the sums $ \sum^J_{j=1, v_{i_j}\in B}v_{i_j}+\nC v_{t+1}$ are all affine lines intersecting $U$. Thus we can certainly have a value $\lambda_{t+1}$ and reset $v_{t+1}=\lambda_{t+1}v_{t+1}$ such that the sum $\sum^J_{j=1, v_{i_j}\in B}v_{i_j}+ v_{t+1}$ are all in $U$ for any $0\leq J\leq k-1$. Then add $v_{t+1}$ into $B$. This finishes the proof.  
	
\end{proof}

\begin{lemma}\label{lm:53} Let $x\in X$ be a generic point. If for a generic vector $v\in T_x(X)$, $\bF_{k,x}(v,\ldots, v)=0$, then the fundamental form $\bF_{k,x}=0$.
\end{lemma}
\begin{proof} Recall the fundamental form $\bF_{k,x}:S^kT_x(X)\rightarrow R^*_{k-1,x}(X)$
	is a symmetric multi-linear map. Let $U\subseteq T_x(X)$ be the open subset in which we choose a generic vector. Then we can have a basis $B=\{v_1,\ldots, v_n\}$ of $T_x(X)$ having the property in Lemma \ref{lm:52} for the number $k$. For any vectors $u_1,\ldots, u_k\in T_x(X)$, $\bF_{k,x}(u_1,\ldots, u_k)$ is a linear combination of the terms of the form $\bF_{k,x}(v_{i_1},\ldots ,v_{i_k})$. Thus it is enough to show $\bF_{k,x}(v_{i_1},\ldots ,v_{i_k})=0$. Write $\widehat{\bF}(v)=\bF_{k,x}(v,\ldots,v)$.	Using the polarization identity (see \cite{Thomas:PolarIdenty}) it suffices to show 
	$$\widehat{\bF}(v_{i_1}+\ldots+v_{i_k})=\widehat{\bF}(v_{i_1}+\ldots+v_{i_{k-1}})=\ldots=\widehat{\bF}(v_{i_1}+v_{i_2})=\widehat{\bF}(v_{i_1})=0 \text{ for all }v_{i_1},\ldots, v_{i_k}\in B.$$
	But this is true by the assumption and the property of $B$.
	
\end{proof}

\begin{proof}[{\bf Proof of Corollary \ref{thm:06}}] We consider the special case of Theorem \ref{thm:vanishing} for $k=1$. It suffices to show that the condition $\mf{t}_1=\theta_1+\delta_1$ is equivalent to $\tau(X)=\Sec(X)$. We calculate the dimension of these two varieties. Note that $\mf{t}_1=n$. By Terracini lemma, the dimension of $\Sec(X)$ is $2n-\delta_1$. One the other hand, let $x$ be a generic point of $X$ and consider the map 
	$$\bar{f}_{2,x}:T_x(X)\times \widehat{\bT}_x(X)/\hat{x}\longrightarrow V^*/\widehat{\bT}_x(X).$$
	By definition, $\theta_1=\rank f_{2,x}(u)=\rank \bar{f}_{2,x}(u)$ where $u\in T_x(X)$ is a generic vector. But we can take an isomorphism $T_x(X)\cong \widehat{\bT}_x(X)/\hat{x}$ so that $\bar{f}_{2,x}(u)=\bII_{x,u}$. Then by Corollary \ref{p:09} the dimension of $\tau(X)$ is $n+\rank \bII_{x,u}$, which equals $2n-\delta_1$. So we conclude $\tau(X)=\Sec(X)$. This proves Corollary \ref{thm:06}. 
	
\end{proof}
\begin{remark} 
	
	The result in Corollary \ref{thm:06} was proved by Griffiths-Harris  \cite[6.15]{Griffiths:AGandLDG} for the case $\dim\tau(X)=2n$ and $r\geq 2n+1$  using the notion of refined third fundamental form. The case of projective nonsingular variety was proved by Landsberg in \cite[10.2]{Landsberg:DegSecTan} using  $II$-generic vectors, refined third fundamental form, and Fulton-Hansen theorem. It is interesting to know if the vanishing of $\bIII_x$ still holds by only assuming $X$ has degenerate secant variety.

\end{remark}

\begin{proof}[{\bf Proof of Corollary \ref{thm:07}}] We prove by contradiction. Assume $\rank \bII_x<\min\{\codim X, \dim X\}$. The second tangent space $\bT_x^2(X)$ has dimension $n+\rank \bII_x$ and thus is a proper linear space in $\nP^r$ of dimension $<2\dim X$. The tangent variety $\Tan(X)=\tau(X)$ is degenerate and thus the secant variety $\Sec(X)$ is also degenerate by Fulton-Hansen theorem. By Corollary \ref{thm:06}, one has the third fundamental form $\bIII_x=0$ and thus $X\subseteq \bT^2_x(X)$, which is a contradiction. 
	
\end{proof}



\bibliographystyle{alpha}

\end{document}